\documentclass[a4paper,fleqn]{cas-dc}

\usepackage[numbers]{natbib}

\usepackage{graphicx}       
\usepackage{subcaption}     
\usepackage{cleveref}

\def\tsc#1{\csdef{#1}{\textsc{\lowercase{#1}}\xspace}}
\tsc{WGM}
\tsc{QE}
\tsc{EP}
\tsc{PMS}
\tsc{BEC}
\tsc{DE}
\newcommand{\bs}{\boldsymbol}

\begin{document}
\let\WriteBookmarks\relax
\def\floatpagepagefraction{1}
\def\textpagefraction{.001}
\shorttitle{SSM Reduction for Geometrically Nonlinear Rotating Structures}
\shortauthors{Gao H.J. et~al.}

\title[mode = title]{Non-intrusive spectral submanifold model reduction for geometrically nonlinear rotating structures with Coriolis and centrifugal forces}                      

\author[1]{Hejun Gao}

\author[2]{Yiliang Wang}

\author[3]{Yan Qing Wang}

\author[4]{Jie Yuan}

\author[1]{Mingwu Li}[orcid=0000-0002-3570-6535]
\cormark[1]
\ead{limw@sustech.edu.cn}

\affiliation[1]{organization={Department of Mechanics and Aerospace Engineering, Southern University of Science and Technology},
                city={Shenzhen},
                postcode={518000}, 
                country={China}}
\affiliation[2]{organization={Independent Researcher},
                city={Pittsburgh},
                country={USA}}                
\affiliation[3]{organization={Key Laboratory of Structural Dynamics of Liaoning Province, College of Sciences, Northeastern University},
                city={Shenyang},
                postcode={110819},
                country={China}}   
\affiliation[4]{organization={Department of Aeronautics and Astronautics, University of Southampton},
                city={Southampton},
                country={UK}}
\cortext[cor1]{Corresponding author}

\begin{abstract}
Rotating structures are widely observed in engineering applications such as turbomachinery and wind turbine. These rotating structures, particularly for blades made by lightweight materials, can undergo large deformation in operations and display complex nonlinear dynamics under the coupling interaction of geometric nonlinearity, Coriolis effect and centrifugal force. Finite element (FE) methods provide a powerful and accurate modeling approach for capturing the complex nonlinear dynamics for realistic rotating structures, yet its high-dimensionality causes significant challenge to efficient prediction for the nonlinear vibration. Here, we present a non-intrusive spectral submanifold (SSM) model reduction for these FE models of rotating structures. We use COMSOL to establish FE models and simulate these FE models to verify the accuracy of SSM-based reduction. We first compute nontrivial static equilibrium induced by the centrifugal force and then construct non-intrusively SSM-based reduced-order model (ROM) anchored at the equilibrium. These SSM-based ROMs enable efficient and accurate extraction of backbone and forced response curves. We use a suite of examples with increasing complexity to demonstrate the effectiveness of the SSM reduction, including a rotating beam, a twisted plate, a rotor with two disks, and an internally resonant fan with three blades. The obtained results also highlight the significant effects of Coriolis force on the nonlinear vibration of rotating structures. 
\end{abstract}




\maketitle

\section{Introduction}
\label{sec:intro} 
Rotating structures are widely used in engineering applications such as aerospace, energy systems, and turbomachinery~\cite{rao1991turbomachine, cassolato1997turbomachinery, rangwala2009structural,xing2025forced,chen2026nonlinear,chai2026response,cui2026dynamic}. Their dynamic behavior plays a critical role in determining the performance, safety, and durability of the overall system. During operation, especially under high rotational speeds, such systems are subjected to multiple interacting physical effects, including centrifugal forces, Coriolis forces, and geometric nonlinearities induced by large deformations~\cite{adams2009rotating, vishwakarma2017vibration, doshi2021review}. The coupling of these effects leads to complex nonlinear dynamic phenomena, such as amplitude-dependent frequency variations, modal interactions, and bifurcations~\cite{genta2005dynamics, luczko2002geometrically, lacarbonara2013nonlinear, LIANG2018393}, making accurate analysis of rotating structures a challenging task.


Accurate analysis of rotating structures is based on mechanics models involving geometric nonlinearity and inertial effects including Coriolis force and centrifugal force~\cite{lai1994nonlinear}. Over the past decades, extensive research has been devoted to the modeling and dynamic analysis of rotating structures~\cite{Ghasabi2018}. Existing approaches can be broadly classified into analytical, semi-analytical, and finite element methods. Analytical and semi-analytical models are particularly valuable for revealing the underlying physical mechanisms of rotating beams, plates and shells~\cite{Ekene2024, Wang2012, santarpia2021hierarchical}. For more complex rotating structures, finite element methods (FEM) provide a more general framework for handling realistic geometries, material distributions and boundary conditions, where all rotational effects, i.e., Coriolis terms, spin-softening effects and geometric stiffening effects are considered~\cite{Azzara2024, Zwlfer2025}.

FEM provides a powerful and accurate modeling method for realistic rotating structures. However, FEM results in high-dimensional, large-scale systems with strong nonlinearities and complex coupling mechanisms. As a consequence, conventional numerical approaches, such as direct time integration, become computationally expensive and inefficient to predict the nonlinear vibration of the rotating structures.

Model order reduction can provide low-dimensional reduced-order models (ROMs) to enable efficient analysis of such high-dimensional systems. Indeed, projection-based reduction methods including modal truncation method~\cite{Liu2002,Beli2018} and proper orthogonal decomposition (POD) method~\cite{Bhartiya2011, Madden2012, Wei2022} have been successfully applied to the model reduction of linear vibration of rotating structures. These linear methods have also been used to the reduction of nonlinear vibration of rotating structures in recent years~\cite{Sakamoto2020, Lu2024}. However, the reduced linear subspaces are not invariant for nonlinear systems~\cite{haller2016nonlinear} and hence the associated ROMs can still be high-dimensional to capture well the nonlinear dynamics. Moreover, POD-based reduction schemes reply on the pre-training data from the full solution which may be very expensive to obtain. To enhance the projection-based method, nonlinear quadratic basis method~\cite{Camacho2025}, hyper-reduction approach~\cite{Kim2022}, and adaptively sampled projections~\cite{Huang2023, Camacho2025} have been applied for nonlinear model reduction of rotating structures as attempts.

Invariant-manifold-based model reduction provides a powerful alternative for constructing ROMs for geometrically nonlinear vibrations~\cite{haller2025modeling,touze2021model}. Linear normal modes span invariant subspaces for the reduction of linear vibration. Under the addition of nonlinear internal force, these invariant subspaces are perturbed into some nonlinear normal modes (NNMs) that are tangent to the subspace~\cite{shaw1991non,shaw1993normal}. Indeed, NNMs have been used to construct ROMs for rotating beams~\cite{pesheck2002modal}. In addition, direct parameterization of NNMs has been recently used to construct ROMs for rotating twisted plate without Coriolis force~\cite{Martin2023}. We note that geometric nonlinearity and centrifugal force are considered in~\cite{pesheck2002modal,Martin2023}, while Coriolis force is not taken into consideration in these studies. The Coriolis force that arises in a non-inertial rotating reference frame can induce gyroscopic coupling and frequency veering. Although Coriolis force may have negligible effects on the linear natural frequencies of rotating structures with certain configuration~\cite{Martin2023}, it plays a significant role in the nonlinear vibration of rotating structures, as we will demonstrate in this study. Therefore, it is essential to account for Coriolis force in the construction of ROMs for rotating structures.

The theory of spectral submanifold (SSM) has been emerged as a rigorous, invariant-manifold-based, exact model reduction framework for high-dimensional nonlinear mechanical systems~\cite{haller2025modeling}. Indeed, there exist infinite many NNMs for a given linear normal mode. Among these NNMs, SSM is defined as the smoothest one under non-resonance conditions~\cite{haller2016nonlinear}. Therefore, SSM is the unique smoothest nonlinear continuation of a given spectral subspace~\cite{haller2016nonlinear}, which provides a natural choice for constructing ROMs thanks to the uniqueness and smoothness. In fact, one can reduces high-dimensional systems to low-dimensional ones with a few degrees of freedom via SSM~\cite{ponsioen2018automated, Ponsioen2020,jain2022compute}, significantly lowering the computational cost of vibration analysis. Recent non-intrusive SSM reduction shows that one can even construct two-dimensional SSM-based ROM for a micro-resonator FE model containing more than a million degrees of freedom~\cite{li2025data}.

Besides finite element based reduction techniques, spectral methods have also attracted increasing attention for nonlinear vibration analysis of rotating structures. In particular, the spectral Chebyshev technique (SCT) has demonstrated excellent convergence properties for rotating beams, plates, and shells while requiring substantially fewer spatial degrees of freedom than conventional finite element discretizations, enabling efficient reduced-order modeling of nonlinear forced vibrations~\cite{Lotfan2023, Mao2025, Lotfan2025}. However, these methods primarily focus on the spatial discretization of continuous systems and require dedicated formulations for different structural configurations. In contrast, the present work aims at a unified, non-intrusive SSM-based reduction framework that applies directly to general finite element models of arbitrary rotating structures within commercial software, achieving model reduction to just a few degrees of freedom without requiring case-by-case derivations.

Here, we apply SSM reduction to geometrically nonlinear rotating structures with Coriolis and centrifugal forces. Indeed, SSM reduction can be applied directly to systems with Coriolis force, as illustrated in the reduction of axially moving beams~\cite{Li2022_P1}, pipes conveying fluids~\cite{li2023nonlinear,wei2026}, and  shaft-disk rotors~\cite{cui2025fully}. Although SSM reduction has been applied to various mechanical systems, it has not been applied to FE models for geometrically nonlinear rotating structures with Coriolis and centrifugal forces. The first goal of this study is to fill this gap by providing such applications. We note that the centrifugal force will result in non-trivial static equilibrium. As the second goal, we will extend the non-intrusive algorithm~\cite{li2025data} to perform SSM reduction at nontrivial equilibrium of rotating structures. Finally, our last goal is to highlight the significance of Coriolis effects on the nonlinear vibration via the SSM-based ROMs.

Both this study and Ref.~\cite{Martin2023} employ invariant manifolds for ROM construction in nonlinear rotating structures, but our work introduces several key novelties relative to~\cite{Martin2023}. First, Ref.~\cite{Martin2023} neglects Coriolis effects, whereas we show that they play a significant role in nonlinear vibrations and cannot be omitted in reduction. Second, their reduction is intrusive, requiring access to explicit multilinear internal-force functions; our approach treats the internal force as a black box, enabling non-intrusive reduction compatible with commercial FE packages like COMSOL. Third, the reduction in~\cite{Martin2023} is limited to simple structures with a single master mode, while we handle complex configurations with internal resonance, including a 1:1:1 resonant rotating fan. Finally, they impose restrictive assumptions on forcing patterns and damping (spatial forcing matching eigenmodes and mass-proportional damping), both of which are lifted in our framework to accommodate arbitrary forcing shapes and general Rayleigh damping.

The rest of this paper is organized as follows. Firstly, the dynamic equation of the rotating model based on the static equilibrium configuration is established through finite element discretization. Then, the SSM method, SSM-based model reduction, and the implementation of its non-intrusive algorithm are introduced. Finally, following the idea from simplicity to complexity, the non-intrusive SSM model reduction is applied to the reduction calculation of various representative rotating structures, demonstrating its computational accuracy and effectiveness. We conclude this paper in \Cref{sec:conclusions}.

\section{Formulation}
\label{sec:formulation} 
The governing equations of rotating structures are formulated using the finite element method (FEM), following a similar approach developed in \cite{GUO2001487}. The model is constructed in a rotating reference frame, and large deformation effects are taken into account. Assuming the rotor rotates around an axis with a constant angular velocity $\hat{\Omega} \in \mathbb{R}^{3}$, the governing equation can be expressed as
\begin{equation}
\textcolor{black}{
\boldsymbol{M}\ddot{\boldsymbol{x}} + (\boldsymbol{C} + \boldsymbol{G}(\hat{\Omega}))\dot{\boldsymbol{x}} + \boldsymbol{F}_{\text{int}}(\boldsymbol{x}) - \boldsymbol{K}_{\text{sp}}(\hat{\Omega})\boldsymbol{x}= \boldsymbol{f}_{\text{cen}}(\hat{\Omega}) + \boldsymbol{f}_{\text{ext}}(t).}
\label{Governing Eqns of Rotor Dyn}
\end{equation}
In Eq.~\eqref{Governing Eqns of Rotor Dyn},
$\boldsymbol{x} \in \mathbb{R}^{n}$ denotes the total displacement vector in the rotating reference frame;
$\boldsymbol{M} \in \mathbb{R}^{n \times n}$ is the mass matrix, $\boldsymbol{C} \in \mathbb{R}^{n \times n}$ is the damping matrix, and $\boldsymbol{G}(\hat{\Omega}) \in \mathbb{R}^{n \times n}$ is the Coriolis matrix. $\boldsymbol{F}_{\text{int}}(\boldsymbol{x}) \in \mathbb{R}^{n}$ represents the internal force vector. $\boldsymbol{K}_{\text{sp}}(\hat{\Omega}) \in \mathbb{R}^{n \times n}$ is the spin softening matrix. $\boldsymbol{f}_{\text{cen}}(\hat{\Omega}) \in \mathbb{R}^n$ denotes the centrifugal force vector evaluated at the reference configuration. $\boldsymbol{f}_{\text{ext}}(t) \in \mathbb{R}^n$ is the excitation force vector. $n$ is the number of degrees of freedom (DOFs) of the global system. We note that $\boldsymbol{f}_{\text{cen}}(\hat{\Omega}) + \boldsymbol{K}_{\text{sp}}(\hat{\Omega})\boldsymbol{x}$ represents the centrifugal force applied on the deformed configuration. For convenience of analysis, we move the term $\boldsymbol{K}_{\text{sp}}(\hat{\Omega})\boldsymbol{x}$ to the LHS of Eq.~\eqref{Governing Eqns of Rotor Dyn}. This treatment is consistent with the standard finite element formulation for rotating coordinate systems implemented in COMSOL~\cite{comsol_rotordynamics_6_2}.

Since the objective of this work is to study vibrations around the static equilibrium state induced by $\boldsymbol{f}_{\text{cen}}$, we decompose the displacement vector $\boldsymbol{x}$ as
\begin{equation}
	\boldsymbol{x} = \boldsymbol{u}_0 + \boldsymbol{u}
	\label{Eqns of TF}
\end{equation}
where $\boldsymbol{u}_0$ denotes the displacement at the static equilibrium state. This displacement satisfies
\begin{equation}
    \boldsymbol{F}_{\text{int}}(\boldsymbol{u}_0) - \boldsymbol{K}_{\text{sp}}(\hat{\Omega})\boldsymbol{u}_0 =\boldsymbol{f}_{\text{cen}}(\hat{\Omega}).
	\label{Eqns of Stat Equil}
\end{equation}
We solve the set of nonlinear algebraic equations~\eqref{Eqns of Stat Equil} using the Newton-Raphson method to obtain the equilibrium $\boldsymbol{u}_0$. In other words, this prestressed equilibrium is obtained through a nonlinear prestressed analysis~\cite{Lotfan2022} under the centrifugal loading. The term $\boldsymbol{u}$ is the dynamic displacement around the static equilibrium state. From now on, we drop the $\hat{\Omega}$ in formulas to make them more compact. By substituting Eq.~\eqref{Eqns of TF} and Eq.~\eqref{Eqns of Stat Equil} into Eq.~\eqref{Governing Eqns of Rotor Dyn}, we obtain
\begin{equation}
	\boldsymbol{M}\ddot{\boldsymbol{u}} + (\boldsymbol{C} + \boldsymbol{G})\dot{\boldsymbol{u}} + \boldsymbol{F}_{\text{int}}(\boldsymbol{u} + \boldsymbol{u}_0) - \boldsymbol{F}_{\text{int}}(\boldsymbol{u}_0) - \boldsymbol{K}_{\text{sp}}\boldsymbol{u} = \boldsymbol{f}_{\text{ext}}(t).
	\label{Eqns of Rotor Dyn - intermediate form}
\end{equation}

To facilitate the application of the SSM computational framework \cite{Li2022_P1,Li2022_P2}, we further decompose the term $\boldsymbol{F}_{\text{int}}(\boldsymbol{u} + \boldsymbol{u}_0) - \boldsymbol{F}_{\text{int}}(\boldsymbol{u}_0)$ into linear and nonlinear components. Let $\boldsymbol{K}_{0}$ denote the tangent stiffness matrix evaluated at the static equilibrium state, which is defined as
\begin{equation}
\boldsymbol{K}_0 = \frac{\partial \boldsymbol{F}_{\text{int}}(\boldsymbol{u})}{\partial \boldsymbol{u}}\bigg|_{\boldsymbol{u}=\boldsymbol{u}_0}.
\label{Eq of Tangent Stiffness Matrix}
\end{equation}
Then the internal force difference can be written as
\begin{equation}
	\boldsymbol{F}_{\text{int}}(\boldsymbol{u} + \boldsymbol{u}_0) - \boldsymbol{F}_{\text{int}}(\boldsymbol{u}_0) = \boldsymbol{K}_{0} \boldsymbol{u} + \boldsymbol{g}_{\text{nl}}(\boldsymbol{u}),
	\label{Taylor Exp of Internal Force}
\end{equation}
where $\boldsymbol{g}_{\text{nl}}$ represents the nonlinear part of the internal force. Let $\boldsymbol{K}_{\text{t}}=\boldsymbol{K}_{0}-\boldsymbol{K}_{\text{sp}}(\hat{\Omega})$ denote the global tangent stiffness matrix, substituting Eq.~\eqref{Taylor Exp of Internal Force} into Eq.~\eqref{Eqns of Rotor Dyn - intermediate form} yields
\begin{equation}
	\boldsymbol{M}\ddot{\boldsymbol{u}} + (\boldsymbol{C} + \boldsymbol{G})\dot{\boldsymbol{u}} + \boldsymbol{K}_{\text{t}}\boldsymbol{u} + \boldsymbol{g}_{\text{nl}}(\boldsymbol{u}) = \boldsymbol{f}_{\text{ext}}(t).
	\label{Eqns of Rotor Dyn - final form}
\end{equation}
In this study, we are interested in the resonance of the vibration and hence the external excitation is assumed to be harmonic. Accordingly,
\begin{equation}
	\boldsymbol{f}_{\text{ext}}(t) = \epsilon \boldsymbol{f}_{\text{load}}\cos{\Omega}t,
	\label{harmonic load}
\end{equation}
where $\Omega$ denotes the excitation frequency and the vector $\boldsymbol{f}_{\text{load}}$ characterizes the spatial distribution of the excitation on the rotating structure. Rayleigh damping is adopted in this study in the form $\bs{C}=\alpha\bs{M}+\beta\bs{K}_\mathrm{t}$, where $\alpha=2\xi_1\omega_1$ and $\beta=0$. We choose $\beta=0$ for consistent comparison with~\cite{Martin2023} but note that SSM reduction works well for general damping~\cite{li2025data}. In each example, $\alpha$ is either prescribed directly or determined from a specified damping ratio $\xi_1$ using the first natural frequency $\omega_1$ defined by the corresponding linear system $\boldsymbol{M}\ddot{\boldsymbol{u}} + \boldsymbol{K}_{\text{t}}\boldsymbol{u} = \boldsymbol{0}$.


\section{SSM-based model reduction}
\label{sec:SSM MOR} 

Now we present how to construct low-dimensional reduced-order models (ROMs) for the full system~\eqref{Eqns of Rotor Dyn - final form} via nonintrusive reduction on spectral submanifolds (SSMs). In particular, we present a brief review on SSM theory, and then provide details on the nonintrusive construction of SSM-based ROMs. We further discuss how to use the SSM-based ROMs to predict backbone and forced response curves. We conclude this section with some implementation details of the SSM-based model reduction.

\subsection{A short review on SSM theory}
We first rewrite~\eqref{Eqns of Rotor Dyn - final form}-\eqref{harmonic load} into the first-order form below
\begin{equation}
\label{eq:1st-ode}
\boldsymbol{B}\dot{\boldsymbol{z}}=\boldsymbol{A}\boldsymbol{z}+\boldsymbol{F}(\boldsymbol{z})+\epsilon\boldsymbol{F}_\mathrm{ext}\cos\Omega t
\end{equation}
where $\boldsymbol{z=(\boldsymbol{x}},\dot{\boldsymbol{x}})$ is the state vector and
\begin{gather}
    \boldsymbol{B}=\begin{pmatrix}\boldsymbol{C} + \boldsymbol{G} & \boldsymbol{M}\\\boldsymbol{M} & \boldsymbol{0}\end{pmatrix}, \boldsymbol{A}=\begin{pmatrix}  -\boldsymbol{K}_{\text{t}} & \boldsymbol{0}\\\boldsymbol{0} & \boldsymbol{M}\end{pmatrix},\nonumber\\
    \boldsymbol{F}=\begin{pmatrix}-\boldsymbol{g}_\mathrm{nl}(\boldsymbol{u})\\\boldsymbol{0}\end{pmatrix},\quad \boldsymbol{F}_\mathrm{ext}=\begin{pmatrix}\boldsymbol{f}_\mathrm{load}\\\boldsymbol{0}\end{pmatrix}.\label{eq:baf}
\end{gather}
We assume that $\boldsymbol{F}(\boldsymbol{z})$ is a smooth nonlinear function such that $\boldsymbol{F}(\boldsymbol{z})\sim\mathcal{O}(|\boldsymbol{z}|^2)$. 

The linear unforced part of~\eqref{eq:1st-ode}, i.e., $\boldsymbol{B}\dot{\boldsymbol{z}}=\boldsymbol{A}\boldsymbol{z}$, admits a trivial equilibrium $\boldsymbol{z}=\boldsymbol{0}$. The associated generalized eigenvalue problem gives
\begin{equation}
    \boldsymbol{A}\boldsymbol{v}_j=\lambda_j\boldsymbol{B}\boldsymbol{v}_j,\quad \boldsymbol{w}_j^\ast\boldsymbol{A}=\lambda_j\boldsymbol{w}_j^\ast\boldsymbol{B},
\end{equation}
where $\lambda_j$ is a generalized eigenvalue and $\boldsymbol{v}_j$ and $\boldsymbol{w}_j$ are the corresponding right and left eigenvectors, respectively. We assume that all eigenvalues have negative real parts and hence the equilibrium is asymptotically stable. We further sort these eigenvalues according to their real parts. With a subset of $M$ eigenvalues with largest real parts, we consider the following $2m$-dimensional master spectral subspace spanned by the associated eigenvectors
\begin{equation}
    \mathcal{E}=\{\boldsymbol{v}_1,\cdots,\boldsymbol{v}_M\}.
\end{equation}
$\mathcal{E}$ is invariant and attracting for the linear system $\boldsymbol{B}\dot{\boldsymbol{z}}=\boldsymbol{A}\boldsymbol{z}$ and one can perform model reduction by projection on $\mathcal{E}$. In practice, we have $M=2$ for systems without internal resonance, and $M=4$ or $M=6$ for systems with two or three normal modes in internal resonance.

The trivial equilibrium persists under the addition of the nonlinear term $\boldsymbol{F}(\boldsymbol{z})$. However, the master subspace $\mathcal{E}$ is not invariant for the nonlinear system. Instead, $\mathcal{E}$ is perturbed into some invariant manifolds that are tangent to $\mathcal{E}$ at the equilibrium. These manifolds are also called as nonlinear normal modes~\cite{shaw1991non}. Under proper non-resonance conditions, it turns out that there exists a unique, smoothest one among these invariant manifolds. This smoothest manifold is defined as the \emph{spectral submanifold} (SSM) associated with $\mathcal{E}$ and denoted by $\mathcal{W}(\mathcal{E})$~\cite{haller2016nonlinear,haller2025modeling}. Since the slow SSM is attracting and invariant, nearby trajectories will converge to $\mathcal{W}(\mathcal{E})$ in exponential rate. Therefore, $\mathcal{W}(\mathcal{E})$ can be used to construct rigorous, low-dimensional, and exact ROMs.

Now we consider the effects of the further addition of the harmonic forcing $\epsilon\boldsymbol{F}_\mathrm{ext}\cos\Omega t$ in~\eqref{eq:1st-ode}. The trivial equilibrium is perturbed into a stable periodic orbit $\gamma_\epsilon$ if $\epsilon$ is small. Accordingly, $\mathcal{W}(\mathcal{E})$ is perturbed as a time-periodic SSM tangent to resonant spectral subbundles of the periodic orbit~\cite{haller2016nonlinear,haller2025modeling}. We denote this time periodic SSM as $\mathcal{W}_{\gamma_\epsilon}(\mathcal{E})$, which can be used to construct ROMs to predict forced responses.

We further discuss the selection of master modes for SSM reduction. For free vibration, the selection is based on the decaying rate and the internal resonance of linear normal modes. In particular, we often take the slowest decaying mode as the master mode because it attracts nearby trajectories. In addition, if the system admits other modes that are internally resonant with the slowest mode, we add these resonant modes into the master modes to capture the internal resonance~\cite{haller2016nonlinear}. As for forced vibration, the selection of master modes is based on both external and internal resonances~\cite{Li2022_P1}. Specifically, we first follow the range of excitation frequency to take the externally resonant mode as the master mode, and then expand the master mode set if there exists other modes that are internally resonant with the excited mode.

We note the reduction on SSM is fundamentally different from modal truncation based reduction. Indeed, SSMs are low-dimensional invariant manifold of the full system~\cite{haller2016nonlinear} and the curved hypersurface of SSMs captures automatically the coupling between slave modes and master modes. In contrast, the modal truncation is to project the full phase space onto lower-dimensional subspace, which neglects the contribution of slave modes and hence can introduce significant errors. These truncation errors, however, are not involved in SSM reduction. Since the reduced dynamics on the SSM often only involves a few degrees of freedom, the SSM reduction can achieve significant speed-up gain for high-dimensional FE problems, as we will demonstrate in \Cref{sec:examples}.


\subsection{Non-intrusive SSM reduction}

In this subsection, we follow~\cite{li2025data} to present how to compute the SSM-based ROMs in a non-intrusive way. More details of this nonintrusive construction can be found in~\cite{li2025data}.

The SSM $\mathcal{W}_{\gamma_\epsilon}(\mathcal{E})$ can be parameterized over the reduced coordinates $\boldsymbol{p}\in\mathbb{C}^M$ as below
\begin{equation}
    \boldsymbol{z}=\boldsymbol{W}_\epsilon(\boldsymbol{p},\phi), \quad \phi=\Omega t
\end{equation}
and its associated reduced dynamics on the SSM is given by
\begin{equation}
    \dot{\boldsymbol{p}}=\boldsymbol{R}_\epsilon(\boldsymbol{p},\phi), \quad \dot{\phi}=\Omega.
\end{equation}
Notably, they satisfy the following invariance equation
\begin{equation}
\label{eq:inv-eq}
\boldsymbol{B}(\partial_{\boldsymbol{p}} \boldsymbol{W}_\epsilon \cdot\boldsymbol{R}_\epsilon+\Omega\partial_\phi \boldsymbol{W}_\epsilon)=\boldsymbol{A}\boldsymbol{W}_\epsilon+\boldsymbol{F}\circ\boldsymbol{W}_\epsilon+\epsilon\boldsymbol{F}_\mathrm{ext}\cos\phi
\end{equation}
where $\circ$ denotes a composition operation for two functions.

To solve for the invariance equation, we first expand the two unknowns in $\epsilon$, yielding
\begin{align}\label{eq:ExpandW_e}
\bs{W}_\epsilon(\bs{p,\bs\phi)} &= \bs{W}\bs{(p)} + \epsilon \bs{X}{(\bs p,\phi)} +\mathcal{O}(\epsilon^2),
\\
\label{eq:ExpandR_e}
 \bs{R}_\epsilon{(\bs p,\bs\phi)} 
&= \bs {R}{(\bs p)} + \epsilon \bs{S} {(\bs p,\phi)} + \mathcal{O}(\epsilon^2).
\end{align}
Here, $\bs{W}\bs{(p)}$ and $\bs {R}{(\bs p)}$ give the manifold parameterization and the vector field of reduced dynamics for \emph{autonomous} SSM at $\epsilon=0$.  We substitute~\eqref{eq:ExpandW_e}-\eqref{eq:ExpandR_e} into the invariance equation~\eqref{eq:inv-eq} and collect $\mathcal{O}(\epsilon^0)$ terms, which gives
\begin{equation}\label{eq:InvEq_aut}
\bs{B} (\text{D} \bs{W}) \bs{R} = \bs{A}\bs{W} + \bs{F} \circ \bs{W}.
\end{equation}
This equation is the invariance equation for the autonomous manifold. 
We further expand $\bs{W}\bs{(p)}$ and $\bs {R}{(\bs p)}$ multivariate monomials, yielding
\begin{equation}
 \label{eq:Rexpansion}
    \bs{W}(\bs{p}) =  \sum_{\bs{m}\in \mathbb{N}^M} \bs{W}_{\bs{m}} \bs{p}^{\bs{m}},
\,\,
    \bs{R}(\bs{p}) =  \sum_{\bs{m}\in \mathbb{N}^M} \bs{R}_{\bs{m}} \bs{p}^{\bs{m}},
\end{equation}
where $\bs{m} \in \mathbb{N}^M$ is a multi-index vector and $\bs{p^m}=p_1^{m_1} \cdots p_M^{m_M}$. In this study, we determine the truncation order of the SSM expansions in~\eqref{eq:Rexpansion} based on the convergence of SSM-based predictions.

We solve for the unknown expansion coefficient vectors $\bs{W}_{\bs{m}}$ and $\bs{R}_{\bs{m}}$ in a recursive fashion. In particular, we substitute the expansions~\eqref{eq:Rexpansion} into~\eqref{eq:InvEq_aut}, and collect terms in the monomial for a given multi-index vector $\boldsymbol{m}$, which gives a system of linear equations for the multi-index vector below~\cite{li2025data}
\begin{equation}
    (\bs{A}-\Lambda_m\bs B)\bs W_{\bs m}=\sum_{j=1}^M \bs{B} \bs{v}_j R^j_{\bs{m}} + \mathcal{C}_\mathbf{m} - \bs{F \circ W}\Bigr|_{\bs{m}}
\end{equation}
where $\Lambda_m=\sum_{i=1}^M m_i\lambda_i$ and the expression of $\mathcal{C}_\mathbf{m}$ is given in~\cite{li2025data}. Here, $ \bs{F \circ W}|_{\bs{m}}$ is a composition polynomial coefficient vector associated with the multi-index vector. Indeed, since $\boldsymbol{F}(\boldsymbol{z})$ and $\boldsymbol{W}(\boldsymbol{p})$ are polynomial functions, their composition is also a polynomial, namely,
\begin{equation}
    \boldsymbol{F}\circ\bs{W}=\bs{F}(\bs{W}(\bs{p}))=\sum_{\bs m}\bs{F \circ W}\Bigr|_{\bs{m}} \bs{p}^{\bs m}.
\end{equation}

In the non-intrusive SSM computation, the coefficient matrices $\bs{B}$ and $\bs{A}$ can be accessed easily because they are based on mass, damping, Coriolis, and stiffness matrices, as seen in~\eqref{eq:baf}. The bottleneck is the non-intrusive evaluation of the composition coefficient vector $ \bs{F \circ W}|_{\bs{m}}$. As detailed in~\cite{li2025data}, this composition coefficient vector can be obtained using only non-intrusive evaluations of the nonlinear function $\bs{F}$, provided that $\bs{F}$ contains only quadratic and cubic nonlinear terms. Indeed, this is the case for this study because the $\boldsymbol{F}_{\text{int}}$ in the equations of motion~\eqref{Governing Eqns of Rotor Dyn} contains only quadratic and cubic terms, which follows from the Green-Lagrange strain used in geometrically nonlinear finite element modeling. Consequently, the $\boldsymbol{g}_{\text{nl}}$ in~\eqref{Taylor Exp of Internal Force} and $\bs{F}$ in~\eqref{eq:baf} also contains only nonlinear terms up to cubic orders.

We follow the same treatment as in~\cite{li2025data} to compute the time-periodic SSM. Specifically, we take a leading-order approximation to the non-autonomous part of the SSM, namely,
\begin{equation}
\label{eq:non-auto-lead}
    \bs{X}(\bs{p},\phi)\approx\bs X_{\bs{0}}(\phi),\quad \bs S(\bs{p},\phi)\approx\bs{S}_{\bs{0}}(\phi).
\end{equation}
We substitute~\eqref{eq:non-auto-lead} into~\eqref{eq:inv-eq} and collect terms at $\mathcal{O}(\epsilon\bs{p}^0)$, yielding
\begin{equation}
\label{eq:leading-order-invariant}
\bs{B}\bs{W}_{\mathbf{I}}\bs{S}_{\bs{0}}(\phi)+\Omega\bs{B}D_{\phi}\bs{X}_{\bs{0}}(\phi)=\bs{A}\bs{X}_{\bs{0}}(\phi)+\bs{F}^{\mathrm{ext}}\cos\phi,
\end{equation}
where $\bs{W}_{\bs{I}}=(\bs{v}_1,\cdots,\bs{v}_M)$. With the ansatz $\bs{X}_{\bs{0}}(\phi)=\bs{x}_{\bs{0}}e^{\mathrm{i}\phi}+\bar{\bs{x}}_{\bs{0}}e^{-\mathrm{i}\phi}$ and $\bs{S}_{\bs{0}}(\phi)=\bs{s}_{\bs{0}}^+e^{\mathrm{i}\phi}+{\bs{s}}_{\bs{0}}^-e^{-\mathrm{i}\phi}$, one can solve for the harmonic coefficients based on the balance of harmonics and also the external resonance condition. More details on solving these harmonic coefficients can be found in~\cite{li2025data}.

\subsection{Prediction via SSM-based ROM}

The SSM-based ROM obtained in the previous subsection can be used to make efficient and even analytic predictions on the free and forced nonlinear vibrations of the full system~\eqref{Governing Eqns of Rotor Dyn} or~\eqref{Eqns of Rotor Dyn - final form}. In particular, we obtain a two-dimensional SSM-based ROM if the system~\eqref{Eqns of Rotor Dyn - final form} admits no internal resonances. Then one can obtain analytic prediction on the backbone and forced response curves of the system via the SSM-based ROM. We refer to~\cite{jain2022compute} for more details on the analytic prediction.

When the linear natural frequencies of the the system~\eqref{Eqns of Rotor Dyn - final form} are internally resonant, the dimensionality of the SSM-based ROM is increased but this dimension only depends on the number of internally resonant modes. For instance, when the first two pairs of modes of the rotating structure admits a near 1:3 internal resonance, we will take these mode as the master subspace $\mathcal{E}$ and then construct a four-dimensional SSM-based ROM. Although analytic predictions are not available in these internally resonant cases, we can make efficient numerical computations since $M\le6$ generally holds yet $n$ can be tens of thousands.

Importantly, the forced periodic orbits of the full system are consistently simplified as fixed points of the SSM-based ROM~\cite{Li2022_P1}. In addition, a rigorous correspondence between the stability and bifurcation of periodic orbits of the full system and that of the fixed points of the SSM-based ROM has been established in~\cite{Li2022_P2}. Therefore, one can infer the stability and bifurcation of the forced periodic orbits easily.


\subsection{Implementation}

\begin{figure}
	\centering
	\includegraphics[width=85mm]{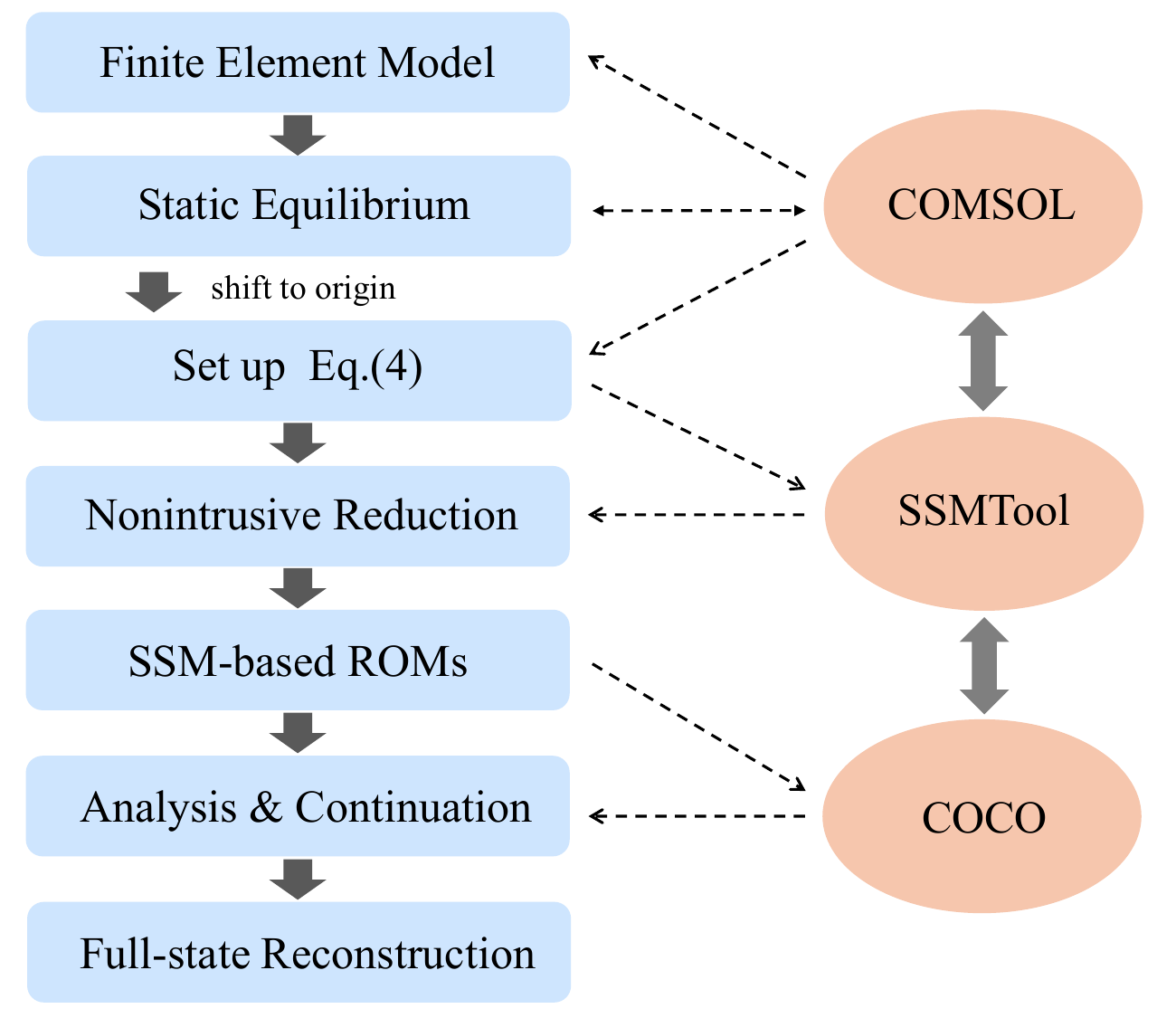}
	\caption{Main steps for the implementation of the nonintrusive SSM-based reduction of rotating structures.}
	\label{fig-flowchart}
\end{figure}

A flowchart that summarizes the main steps of our implementation is given in \Cref{fig-flowchart}.
We use COMSOL Multiphysics to perform finite element modeling of rotating structures in this study. With a given constant angular velocity $\hat{\Omega}$, we first solve for the static equilibrium $\boldsymbol{u}_0$ induced by the centrifugal force in COMSOL. In practice, one can apply Newton-Raphson scheme to Eq.~\eqref{Eqns of Stat Equil} to extract the static equilibrium.

With the static equilibrium shifted to origin, we then set up the full system~\eqref{Eqns of Rotor Dyn - final form} as a pre-processing for the SSM reduction. Specifically, COMSOL allows for the extraction of the mass matrix $\bs{M}$, effective damping matrix $\bs{C}+\bs{G}$, the tangent stiffness matrix $\bs{K}_t$ and also the internal force function $\bs{F}_\mathrm{int}$~\cite{comsol-matlab}. We then use the relation in~\eqref{Taylor Exp of Internal Force} to construct a function handle for $\bs{g}_\mathrm{nl}$.

We further use SSMTool~\cite{SSMTool2} to perform the non-intrusive computation of the SSM and the associated ROM. As documented in~\cite{li2025data}, the non-intrusive computation has been implemented in the model reduction package SSMTool, which supports automated computation of arbitrary dimensional SSMs up to any order of truncation. Indeed, we need to truncate the Taylor series~\eqref{eq:Rexpansion} in practice. In this study, we perform such truncation based on the check of convergence for these Taylor series. An alternative approach to truncate is to use the residual of invariance equation~\cite{li2023model}. More details on the coupling of SSMTool and COMSOL can be found in~\cite{li2025data}.

A numerical continuation package, \textsc{coco}~\cite{Dankowicz2013RecipesContinuation,Ahsan2022MethodsEquations}, has been integrated into SSMTool~\cite{Li2022_P2}. This integration enables effective stablity and bifurcation analysis of the fixed points of the SSM-based ROM, which further enables us to infer the stabilty and bifurcation of the forced periodic orbits of the full system. Here, we use the parameter continuation in \textsc{coco} to extract the forced response curve of periodic orbits and detect bifurcations on them.

\section{Examples}\label{sec:examples} 
In this section, the proposed SSM-based model reduction is applied to three representative rotating structures. Backbone curves and FRCs are computed to assess the accuracy of the method. The results are compared against reference solutions from the literature when available; otherwise, direct time integration of Eq.~\eqref{Eqns of Rotor Dyn - final form} is performed to obtain steady-state responses as benchmark solutions. Most models are discretized using 27-node hexahedral elements, except for \Cref{sec:fan} which uses 10-node tetrahedral elements. COMSOL is also used for direct time integration of Eq.~\eqref{Eqns of Rotor Dyn - final form} to obtain the system response.

\subsection{Rotating Beam}
\label{sec:beam} 
In this example, we perform SSM-based reduction for a rotating beam to investigate different nonlinear dynamic scenarios. First,~\Cref{sec: without Coriolis effect} considers a baseline case without the Coriolis effect to reproduce and validate results available in the literature. Next,~\Cref{sec: considering Coriolis effect} analyzes the complete model including the Coriolis force to assess its influence on the system dynamics. Finally, in~\Cref{sec: beam inner resonance} the rotational speed is tuned to induce internal resonance, and our method is applied to evaluate its performance in more complex nonlinear regimes.

\begin{figure}
	\centering
	\includegraphics[width=85mm]{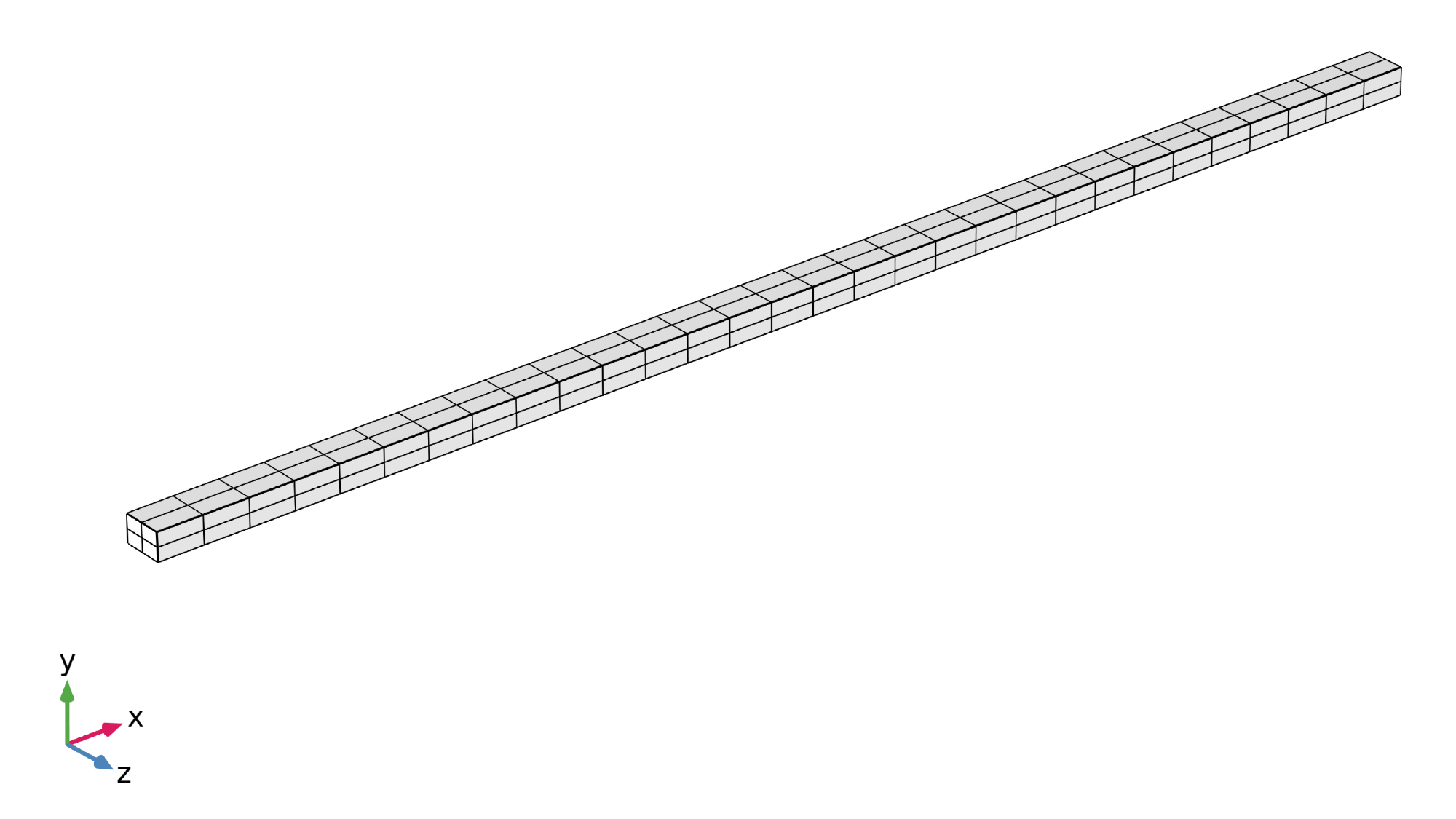}
	\caption{Geometry of the rotating beam model and its finite element mesh.}
	\label{fig-beam}
\end{figure}

The geometry and corresponding finite element mesh of the rotating beam are adopted from \cite{Martin2023}. \Cref{fig-beam} illustrates the geometry and mesh of the rotating cantilever beam model, with its left end face clamped. The beam has a length of 1 m (in the x direction), a thickness of 0.02 m (in the y direction), and a width of 0.03 m (in the z direction). The z-axis is defined as the rotation axis, located 0.1 m away from the clamped end of the beam. The beam material has a Young's modulus of 104 $\mathrm{GPa}$, a Poisson's ratio of 0.3 and a mass density of 4400 $\mathrm{kg/m^3}$. The degrees of freedom (DOFs) of this model is 4500.

\subsubsection{Model without Coriolis effect}\label{sec: without Coriolis effect} 
In~\cite{Martin2023}, the direct parameterize method of invariant manifolds (DPIM) was developed and applied to analyze the rotating beam model. To focus primarily on the centrifugal effect, the Coriolis effect is not considered. For consistency, we consider the same configuration by neglecting the Coriolis term $\boldsymbol{G}(\hat{\Omega})$ and damping term in Eq.~\eqref{Eqns of Rotor Dyn - final form}. We take the first bending mode as the master subspace to construct two-dimensional SSM-based ROM. We truncate the SSM expansion at $\mathcal{O}(7)$, consistent with the expansion order used in~\cite{Martin2023} for DPIM.

\begin{figure}
	\centering
	\includegraphics[width=85mm]{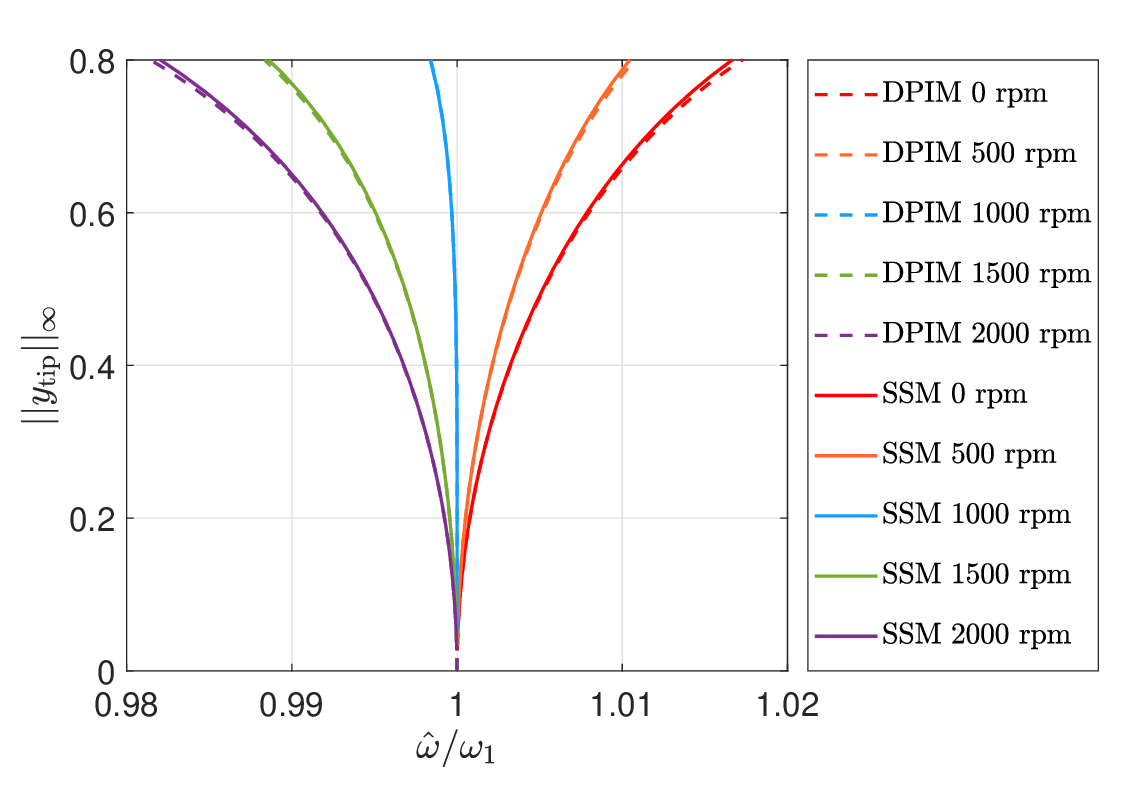}
	\caption{Backbone curves of the rotating beam obtained via different methods with Coriolis effect neglected. Here, $\hat{\omega}$ represents the response frequency, $\omega_1$ denotes the frequency of first bending mode, and the vertical axis indicates the deflection at the beam tip.}
	\label{fig-beam-bb-ref}
\end{figure}

With the SSM-based reduction, backbone curves are computed for a range of rotational speeds (0, 500, 1000, 1500, and 2000 rpm), as shown in \Cref{fig-beam-bb-ref}. The SSM results are in agreement with those reported in \cite{Martin2023}, demonstrating the accuracy of our proposed approach when the Coriolis effect is neglected. In addition, the characteristic hardening–softening transition behavior of the rotating beam is captured \cite{Martin2023, Thomas2016, luo2025role, Zhao2017, Bera2023, Turhan2009}.

\subsubsection{Model considering Coriolis effect}\label{sec: considering Coriolis effect} 
Then, the rotating beam is modeled without neglecting the Coriolis force term in Eq.~\eqref{Eqns of Rotor Dyn - final form}. The proposed SSM based ROM is applied to this model to obtain the backbone curves and FRCs.

First, three rotational speeds (500, 1000, and 2000 rpm) are selected to compute the backbone curves of the model with the Coriolis force included. The damping term is neglected here, as only the autonomous reduced-order model is considered for computing the backbone curves. The results and their comparison with those obtained in \Cref{sec: without Coriolis effect}, where the Coriolis force was neglected, are shown in \Cref{fig-beam-bb-Coriolis}. The Coriolis force has a limited influence on the natural frequencies of the system (corresponding to the zero-amplitude points on the backbone curves), but it significantly affects the nonlinear dynamic response. {\Cref{tab-beam} details the comparison of vibration amplitudes on the backbone curves at sampled dimensionless frequencies. As seen in the table, the relative error of neglecting Coriolis force is large and can even exceed 145\% at $\Omega=1000$ rpm. In general, the Coriolis force induces obvious differences in the nonlinear dynamic responses of the model, and such discrepancies increase as the vibration amplitude grows. Therefore, when performing dynamic response analysis for such rotating structures, the inclusion of the Coriolis force should be carefully considered.

\begin{figure}
	\centering
	\includegraphics[width=85mm]{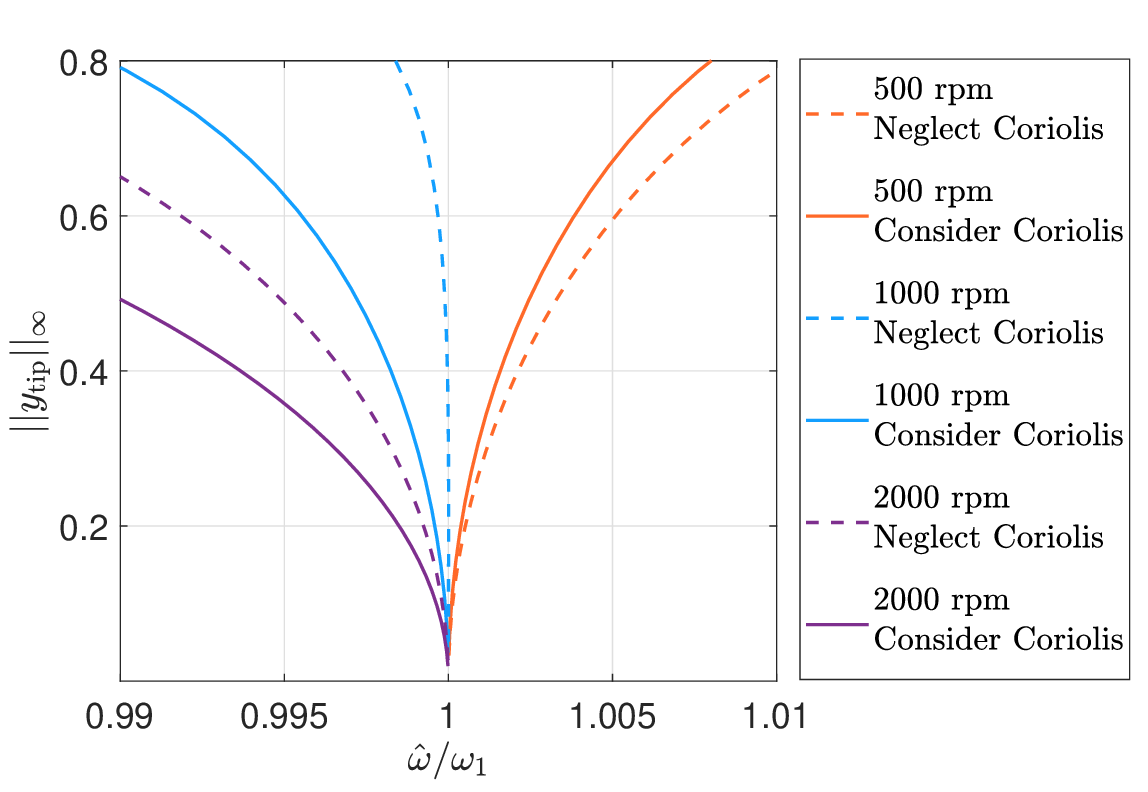}
	\caption{Backbone curves of the rotating beam obtained via SSM reduction. The meanings of  $\hat{\omega}$, $\omega_1$ and the vertical axis are the same as those of \Cref{fig-beam-bb-ref}.}
	\label{fig-beam-bb-Coriolis}
\end{figure}

\begin{table*}[cols=5,pos=h]
\caption{Comparison of vibration amplitudes on the backbone curves of the beam model (see Fig.~\ref{fig-beam-bb-Coriolis}) at dimensionless sampled frequencies $\hat{\omega}/\omega_1$ with and without Coriolis force. Here, the relative error takes the amplitude with Coriolis force (Amp. w/Coriolis) as the reference solution and compute the error for the amplitude without Coriolis force (Amp. w/o Coriolis).}
\label{tab-beam}
\begin{tabular*}{1\textwidth}{@{\extracolsep{\fill}} ccccc @{}}
    \toprule
    Rotational Speed (rpm) & $\hat{\omega}/\omega_1$ & Amp. w/ Coriolis (m) & Amp. w/o Coriolis (m) & Relative Error\\
    \midrule
    500  & 1.005 & 0.663 & 0.608 & 8.2\%\\
    1000 & 0.999 & 0.294 & 0.722 & 145.6\%\\
    1500 & 0.995 & 0.365 & 0.486 & 33.2\%\\
    \bottomrule
\end{tabular*}
\end{table*}

Furthermore, a converged FRC of rotating beam model at a high rotational speed of 2000 rpm is obtained via non-intrusive SSM reduction up to $\mathcal{O}(5)$ truncation. Since no reference results are available in the literature, direct time integration of \eqref{Eqns of Rotor Dyn - final form} is performed to obtain the dynamic response of the full model, which serves as the benchmark solution. Damping parameter $\alpha$ is set to be 0.5, and a harmonic point load of the form $26\cos\Omega t$ is applied at the center of the beam end such that the maximum transverse displacement reaches approximately 50\% of the beam length, thereby inducing strong geometric nonlinearity. For comparison, results obtained without including the Coriolis force are also presented. For the reference solutions, a set of discrete excitation frequencies is selected, and the corresponding steady-state responses are computed individually. All results are shown in \Cref{fig-beam-FRC}.

\begin{figure*} 
    \centering
	\includegraphics[width=170mm]{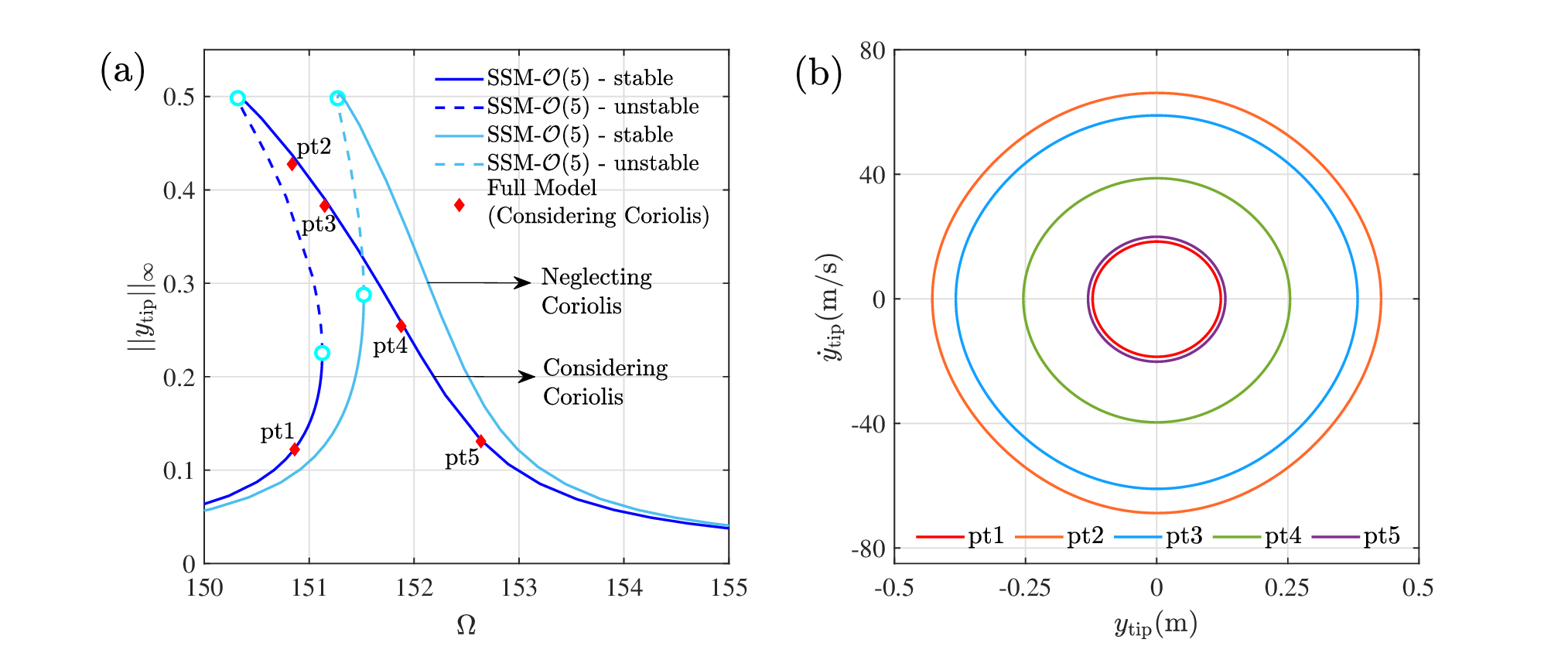}
	\caption{The FRCs of rotating beam model obtained via SSM model order reduction, along with the comparison of Coriolis force considering and Coriolis force neglecting, as well the results of full model. The left picture (a) shows FRC of the rotating beam under 2000 rpm, with five representative data points of the full model as the red marks. Here and throughout this paper, the cyan circles denote saddle–node bifurcation points. The right picture (b) is the phase portrait corresponding to the five full model data points in the left picture, highlighting that each point on the FRC corresponds to a periodic response of the full model.}
	\label{fig-beam-FRC}
\end{figure*}

We observe from \Cref{fig-beam-FRC} the significant effects of Coriolis force on the FRC, especially at large vibration amplitudes. This observation is consistent with the case of backbone curves. Notably, the resonance peak shifts only slightly when the Coriolis force is neglected: from $\Omega_{\mathrm{res,w}}=150.320$ to $\Omega_{\mathrm{res,w/o}}=151.271$, a change of about 0.6\%. Despite this seemingly small frequency shift, the corresponding change in vibration amplitude at the same frequency is dramatic. For example, at $\Omega=\Omega_{\mathrm{res,w}}$, the model neglecting Coriolis force predicts an amplitude of 0.07, whereas the model including the force gives 0.50-a difference of more than 7 times. This large discrepancy at nearly the same frequency highlights that the Coriolis force can strongly affect the vibration amplitude near resonance.

As shown in~\Cref{fig-beam-FRC}, the reduced order model solutions agree well with the full model solutions obtained via direct time integration. Notably, obtain only one steady-state response in \Cref{fig-beam-FRC}(a) takes more than 3 hours, whereas generating a complete FRC using our proposed reduced order model requires less than 17 minutes in total. This demonstrates that model order reduction significantly reduces computational cost while maintaining accuracy. In fact, we have used initial conditions from the SSM-based prediction when performing the time integration of the full system, which saves a significant amount of time to reach the steady state for the simulation. We use this scheme for choosing initial conditions for direct time integration of full systems throughout this study.

\subsubsection{Rotating beam with inner resonance}\label{sec: beam inner resonance}
To show that even for rotating structure with internal resonance, SSM can still effectively perform model reduction on the system and obtain accurate results, this section investigates the nonlinear vibration characteristics of rotating beam with 1:3 internal resonance at a specific rotational speed. 

Due to the centrifugal effect, the frequencies of most modes of the rotating beam model increase with the rise of rotational speed, and the rising rates are not consistent \cite{Martin2023, Wright1982}. This makes it possible for internal resonance to occur between different modes at a specific rotational speed. Here, we add an isotropic spring attached to the free end of the beam to trigger the internal resonance easily. This modification significantly increased the frequency of first bending mode. When the spring stiffness is $33,798\,\text{N/m}$ and the rotational speed is $100\,\text{rad/s}$, a 1:3 multiple relationship is formed between the first and third bending modes, which indicates that the system may generate 1:3 internal resonance when nonlinearity is considered. In this case, damping ratio of first mode $\xi_1$ is assigned as 0.001.

\begin{figure*}
	\centering
	\includegraphics[width=170mm]{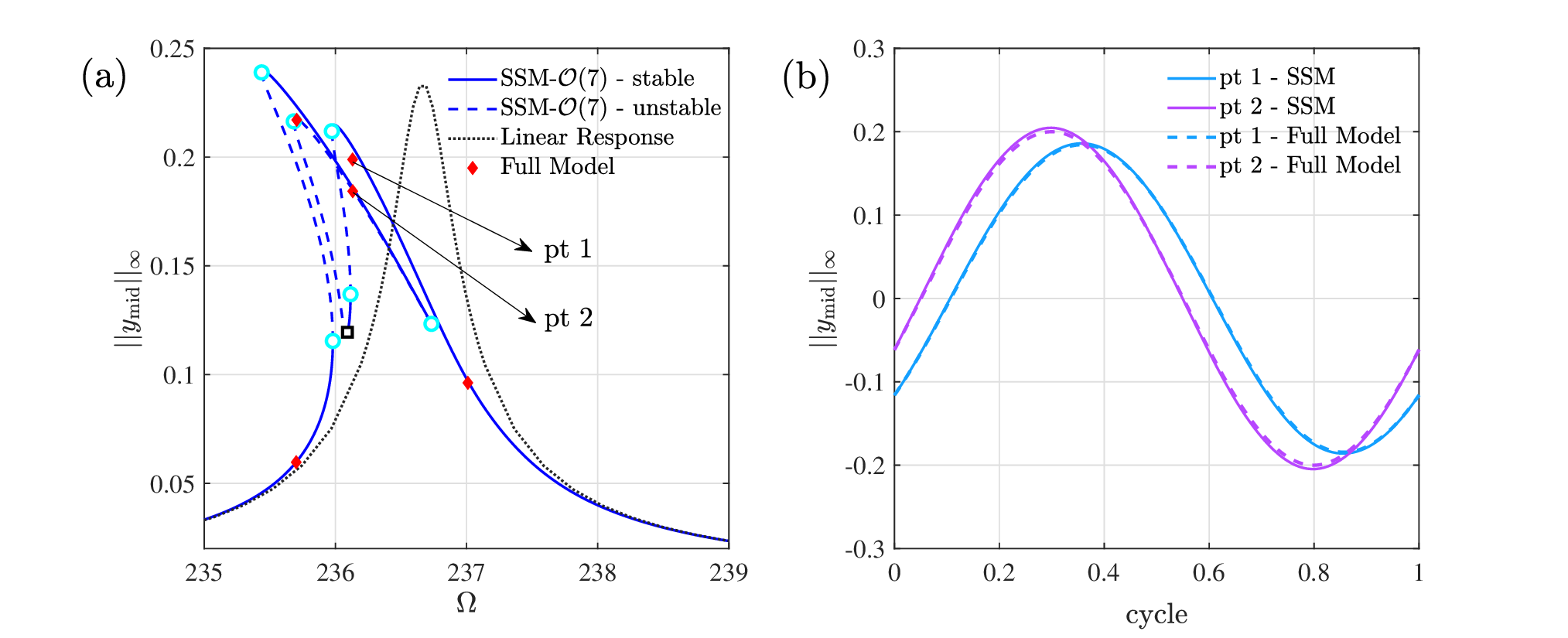}
	\caption{The forced response curve (FRC) of rotating beam with inner resonance obtained via SSM model order reduction, along with the results of full model. The left panel (a) shows FRC of the rotating beam under 1:3 internal resonance, with five representative data points of the full model as the red marks. Specially, here the black square denotes a Hopf bifurcation point. The right panel (b) is the time response series of the two picked points in the left picture, with comparison of SSM reduction and full model.}
	\label{fig-beam-resonance}
\end{figure*}

We take the internally resonant modes as the master subspace and construct the associated four-dimensional SSM-based ROM. \Cref{fig-beam-resonance} shows the FRC of the system obtained from the SSM-based ROM. Meanwhile, several representative excitation frequencies are selected to obtain the responses of full model, which are indicated by the red marks in the left picture of \Cref{fig-beam-resonance}.
The results show that the reduced order model is in good agreement with full model, with a relative error not exceeding 3\%. At the same time, internal resonance makes the response behavior of the nonlinear system more complex, which is inherently different from that of reduction on SSM with single master mode (cf.~\Cref{fig-beam-FRC}). For example, within certain ranges of excitation frequencies, two sets of stable solutions with similar amplitudes appear, and there is a phase difference between their responses, as shown in the right panel of \Cref{fig-beam-resonance}.

\subsection{Rotating Twist Blade}
\label{sec:blade} 

The rotating beam model in \Cref{sec:beam} has structural symmetry, and only odd-order nonlinear terms arise when the nonlinear terms in the governing equations are expanded as a power series \cite{Ekene2024}. To further verify the applicability of the SSM method to complex structures, a twist blade model is considered here. In this model, the structural symmetry is broken, leading to a more general form of the nonlinear terms in the governing equations \cite{sinha2011natural, Yao2011-1, Li2023-1}.  Damping ratio of first mode $\xi_1$ is assigned as 0.005.

\Cref{fig-blade} illustrates the geometry and finite element mesh of twist blade model. The total DOFs of this model is 25500. The model has a rectangular cross-section and a length of $1\,\text{m}$, with a total torsion angle of $60^\circ$ between the root and the tip, which is linearly distributed along the length. The rotation axis is $0.3\,\text{m}$ from the root and its direction is along $(0,0,1)$ and the left surface at $x=0.3\,\text{m}$ is clamped and the other three surfaces are free. The material properties are the same with those of the rotating beam model described in \Cref{sec:beam}.

\begin{figure}
	\centering
	\includegraphics[width=85mm]{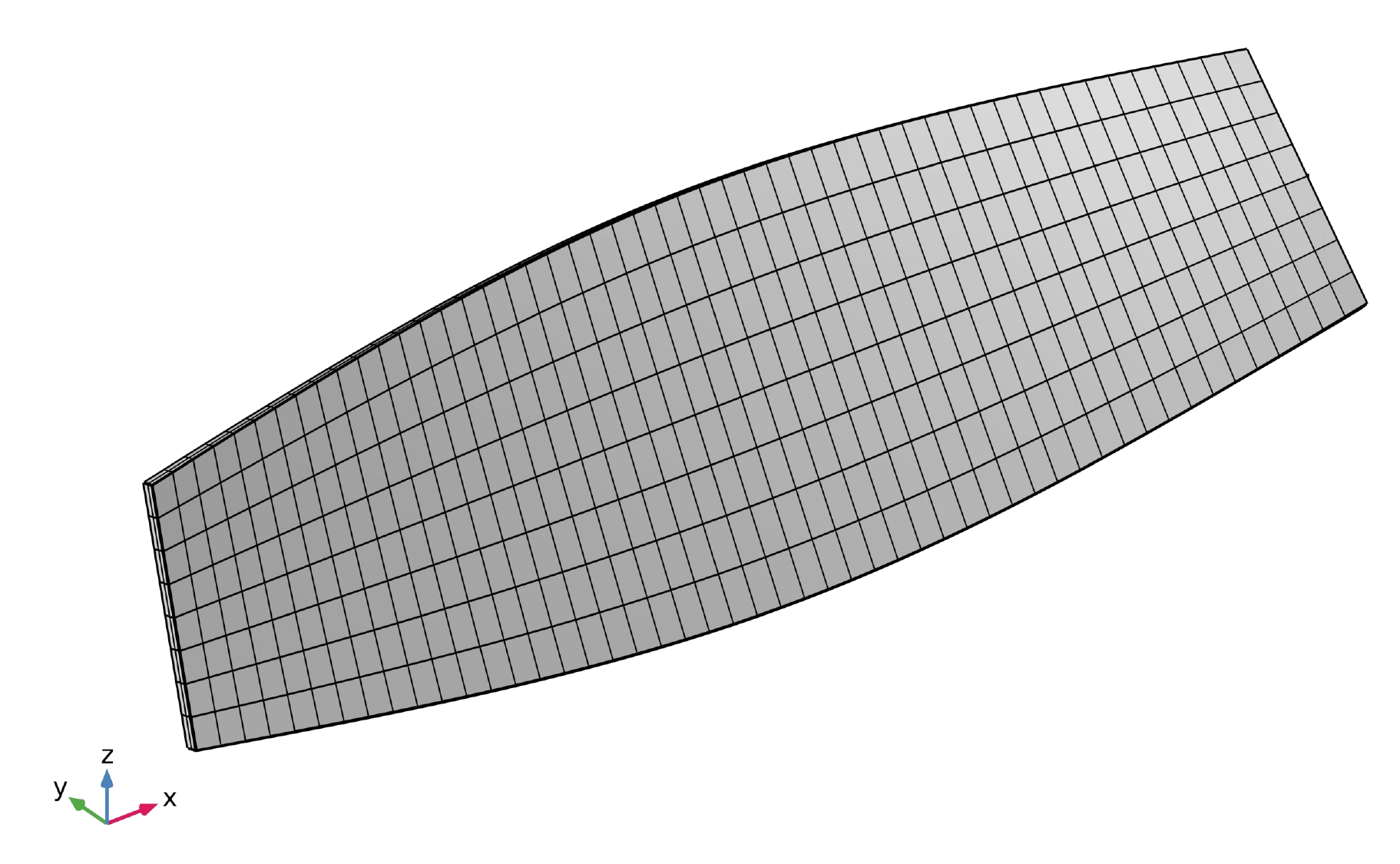}
	\caption{Geometry of the twist blade model and its finite element mesh.}
	\label{fig-blade}
\end{figure}

Here, we take the first mode as the master subspace to perform model reduction on the associated two-dimensional SSM and then use the SSM-based ROM to study the primary resonance for the nonlinear vibration of blade model. In this case, an $\mathcal{O}(5)$ truncation of SSM expansion is sufficient to yield converged predictions. Based on the SSM-based ROM, the backbone curves of the system are computed with and without considering the Coriolis force at different rotational speeds (300, 500, and 1000 rpm), as shown in \Cref{fig-blade-bb}. 

\begin{figure}
	\centering
	\includegraphics[width=85mm]{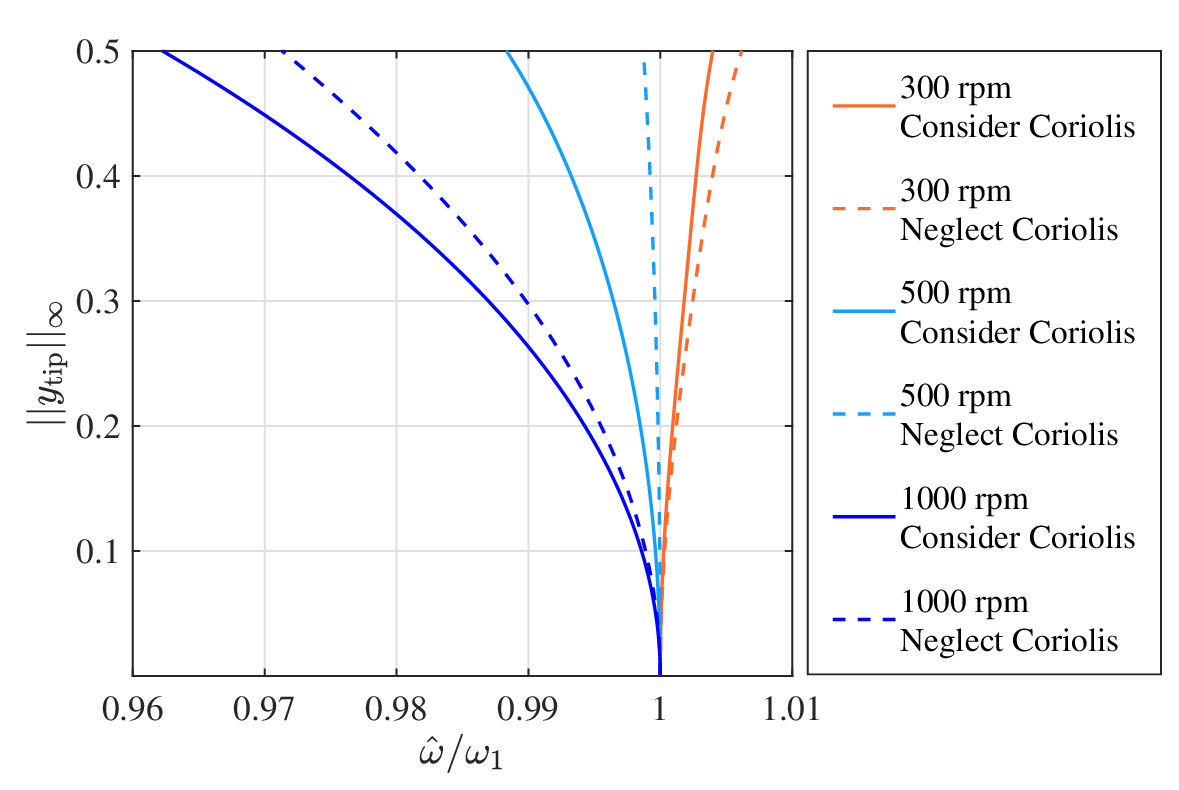}
	\caption{Backbone curve of the twist blade under different rotation speeds.}
	\label{fig-blade-bb}
\end{figure}

Similar to the case of rotating beam, we infer from the backbone curves that the twisted blade structure also exhibits a transition of the system from nonlinear hardening to nonlinear softening as the rotational speed increases. In addition, a significant difference is consistently observed between the backbone curve with and without the Coriolis force when the excitation frequency $\Omega$ deviates from the linear natural frequency, further highlighting the importance of including the Coriolis effect.

We further study the forced vibration of the twist plate at a relatively high rotational speed (1000 rpm). A harmonic point load of the form $150\cos\Omega t$ is applied at the center of the plate end. We obtain the FRC of the system via the SSM-based ROM and present the results in \Cref{fig-blade-FRC}. We note that the peak vibration amplitude on the FRC exceeds 40\% of the plate length, indicating the presence of strong geometric nonlinearity. Indeed, the presence of unstable periodic orbits on the FRC indicates a strongly nonlinear system response.

We again select several representative excitation frequencies and perform direct numerical integration of the full system to validate the accuracy of the SSM-based prediction. As shown in \Cref{fig-blade-FRC}, even at large amplitudes, the reduced-order model is in excellent agreement with the reference solutions obtained via direct time integration, with a relative error not exceeding 2\%. Here, the computational time is around $5\,\text{mins}$ for $\mathcal{O}(3)$, $38\,\text{mins}$ for $\mathcal{O}(5)$, whereas obtaining only one steady-state response takes around 5 hours, which again demonstrates the significant speed-up gain of the FRC.

\begin{figure}
	\centering
	\includegraphics[width=85mm]{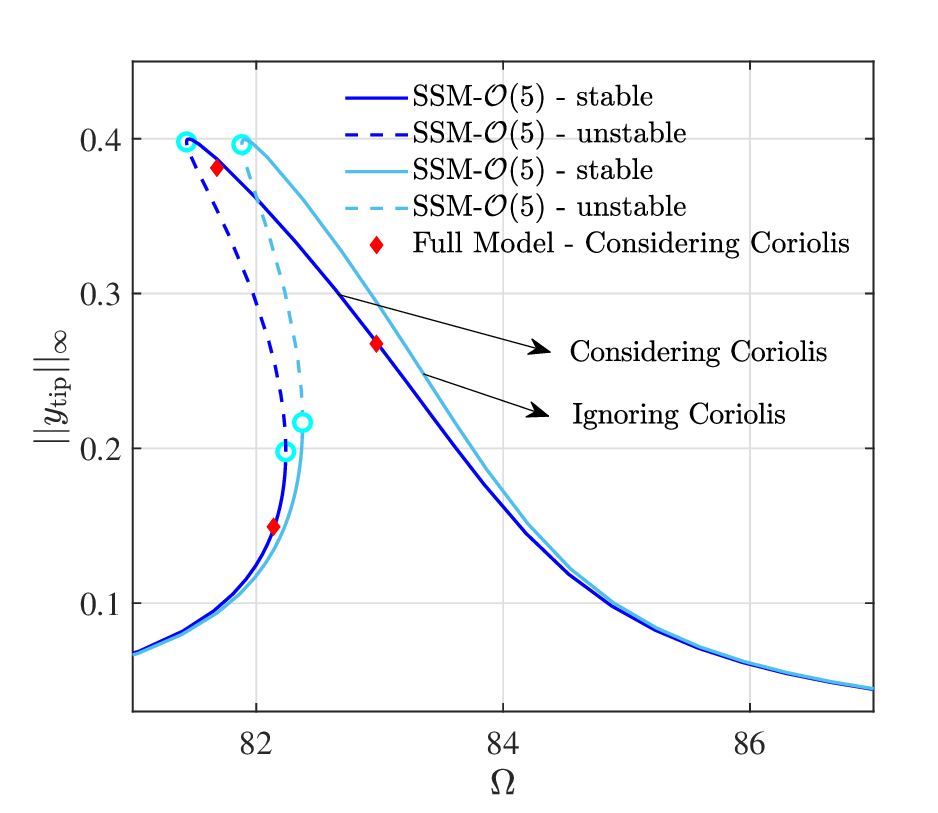}
	\caption{FRC of twist blade under 1000 rpm.}
	\label{fig-blade-FRC}
\end{figure}

\subsection{Rotor with Two Disks}
\label{sec:rotor sys} 

We also apply the proposed SSM-based ROM to a rotating structure with multiple components. The model consists of a flexible shaft and two flexible disks. In this example, the nonlinear dynamic response of a simple rotor system is investigated at a relatively low rotational speed, where the rotor is not yet buckled.

Consider the model shown in \Cref{fig-rotor}: the shaft has a length of $1\,\mathrm{m}$ and a radius of $1\,\mathrm{cm}$. Two metal disks with radio of $8\,\mathrm{cm}$ and $6\,\mathrm{cm}$, both with a thickness of $2\,\mathrm{cm}$, are mounted at distances of 0.3 m and 0.7 m from one end of the shaft, respectively. The larger disk is laterally offset by $5\,\mathrm{mm}$ in the z-direction. The shaft is made of steel, with a density of $7850\,\mathrm{kg/m^3}$, a Young’s modulus of $210\,\mathrm{GPa}$, and Poisson’s ratio of $0.3$. The disks are made of aluminum with a density of $2700\,\mathrm{kg/m^3}$, a Young’s modulus of $72\,\mathrm{GPa}$, and a Poisson’s ratio of $0.3$. Both ends of the shaft are clamped, and the rotational speed is $100 \,\mathrm{rad/s}$. The geometry and material properties of this model are adopted from~\cite{Holzinger2025}. Damping ratio of first mode $\xi_1$ is assigned as 0.005. The total DOFs of this model is 47673.

\begin{figure}
	\centering
	\includegraphics[width=85mm]{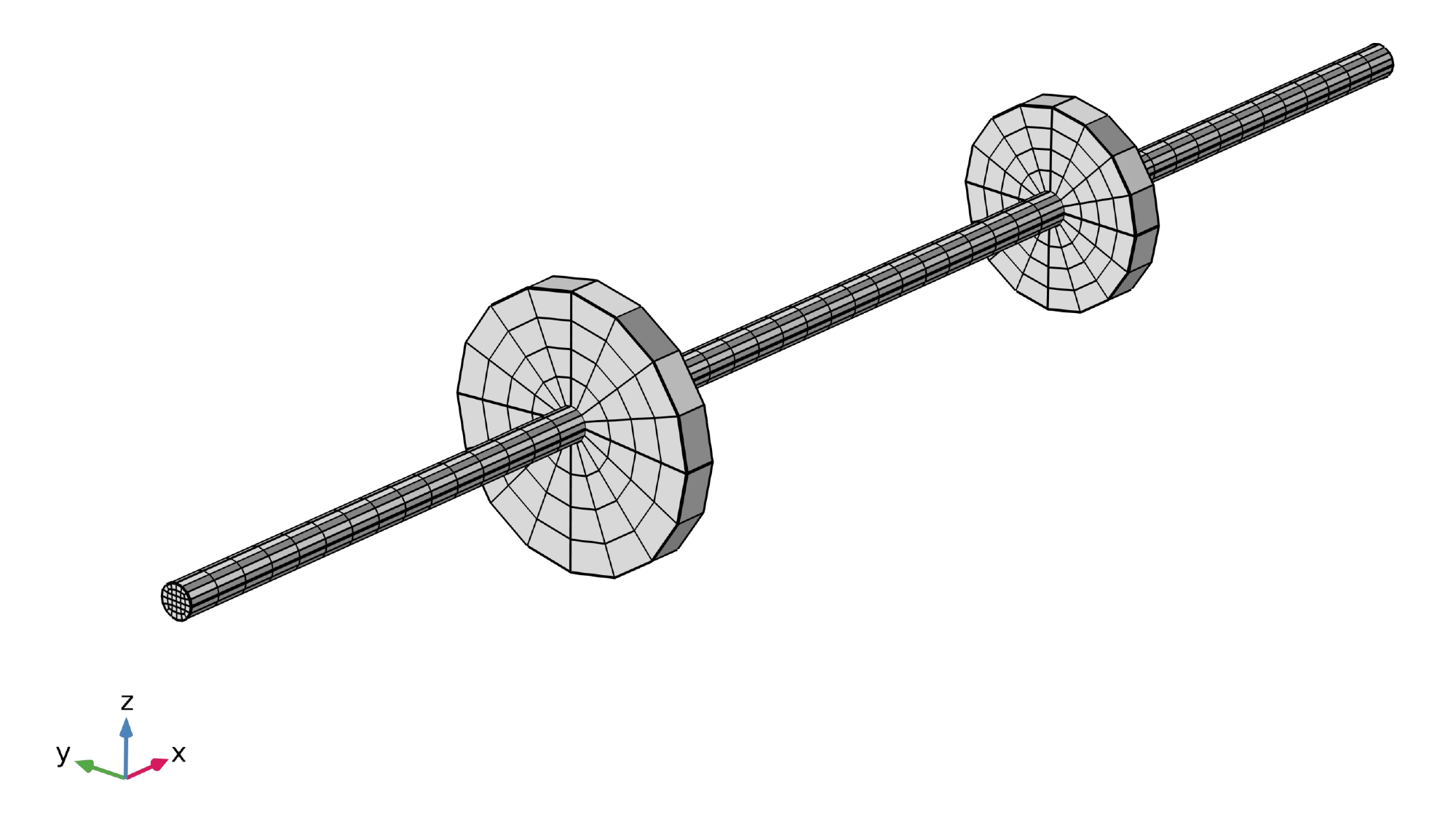}
	\caption{Geometry of the rotor model with two disks and its finite element mesh.}
	\label{fig-rotor}
\end{figure}

\begin{figure}
	\centering
	\includegraphics[width=85mm]{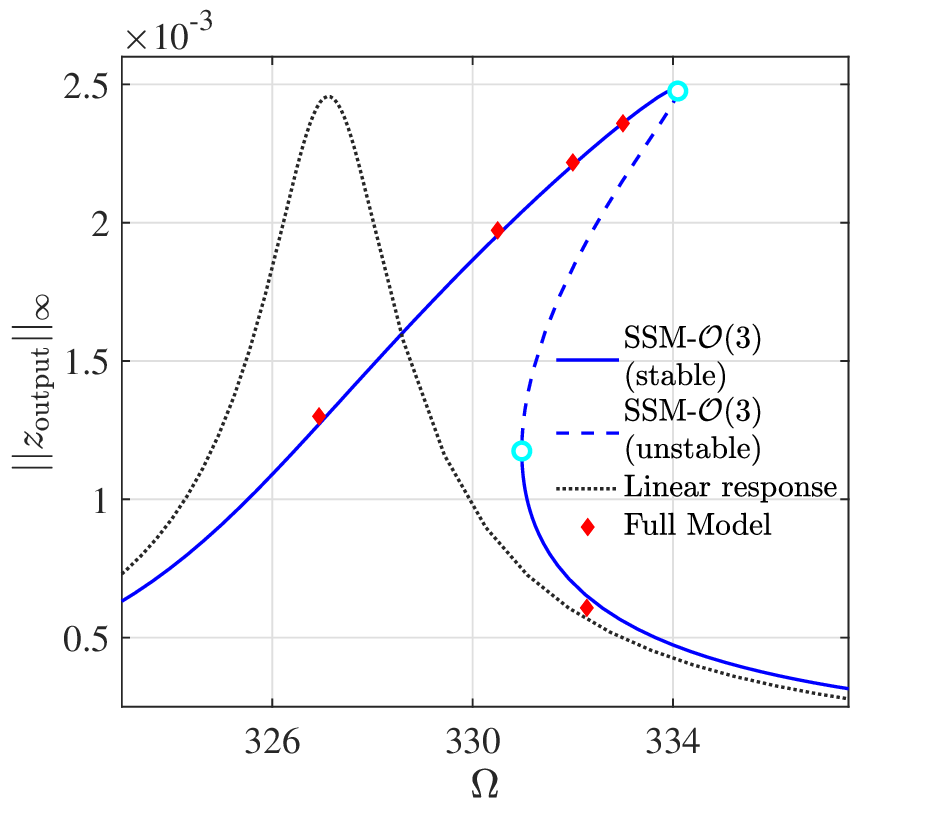}
	\caption{FRC of simple rotor model under 100 rad/s.}
	\label{fig-rotor-FRC}
\end{figure}

\begin{table}[width=1\linewidth,cols=4,pos=h]
\caption{The first three natural frequencies in rad/s of the rotor model with or without consideration of Coriolis force.} 
\label{tab-rotor}
\begin{tabular*}{\tblwidth}{@{}LLLL@{} }
    \toprule
    & Mode 1 & Mode 2 & Mode 3 \\
    \midrule
    w/ Coriolis force & 327.13 & 525.32 & 940.31 \\
    w/o Coriolis force & 414.52 & 414.58 & 1033.96 \\
    \bottomrule
\end{tabular*}
\end{table}

Unlike the previous two examples, where the Coriolis force has minor effects on the linear natural frequencies of the rotating structures~\cite{Martin2023}. Here, whether the Coriolis force is considered or not can lead to considerable discrepancies in the natural frequencies of the model, as shown in \Cref{tab-rotor}. Indeed, the $z$-motion for the disks results in Coriolis force along the $y$-direction, which further affects the $y$-motion of the disks. Therefore, it is essential to consider the Coriolis force even at linear dynamics. In the rest of this example, we take the Coriolis force into consideration.

To investigate the forced vibration of the system, harmonic loads along the z‑axis direction are applied at (0.3, 0, 0.075) and (0.7, 0, 0.06), corresponding to the edges of the two disks. Since no internal resonances are observed, we take the first mode as the master subspace to perform two-dimensional SSM reduction. In particular, a non-intrusive SSM reduction up to $\mathcal{O}(3)$ truncation is performed to yield a converged FRC presented in \Cref{fig-rotor-FRC}. Here, the $z$-displacement at the locations where the first load applied is chosen as the output for the FRC. We observe that the FRC displays a hardening behavior and it differs significantly from the linear response near the resonance region. 

Similar to the previous examples, we use direct time integration of the full system to validate the accuracy of the SSM-based prediction. Several representative excitation frequencies are selected, and the corresponding steady-state responses obtained via the direct time integration responses are shown as red markers in \Cref{fig-rotor-FRC}. The results show good agreement, which validates the accuracy of the SSM predictions. We again observe a significant speed-up gain from the SSM-based reduction. Specifically, computing a single data in \Cref{fig-rotor-FRC} takes more than 8 hours, whereas generating a complete FRC using our proposed reduced order model requires less than 7 minutes in total.

\subsection{Rotating Fan with 1:1:1 Internal Resonance}
\label{sec:fan} 

\begin{figure}
	\centering
	\includegraphics[width=85mm]{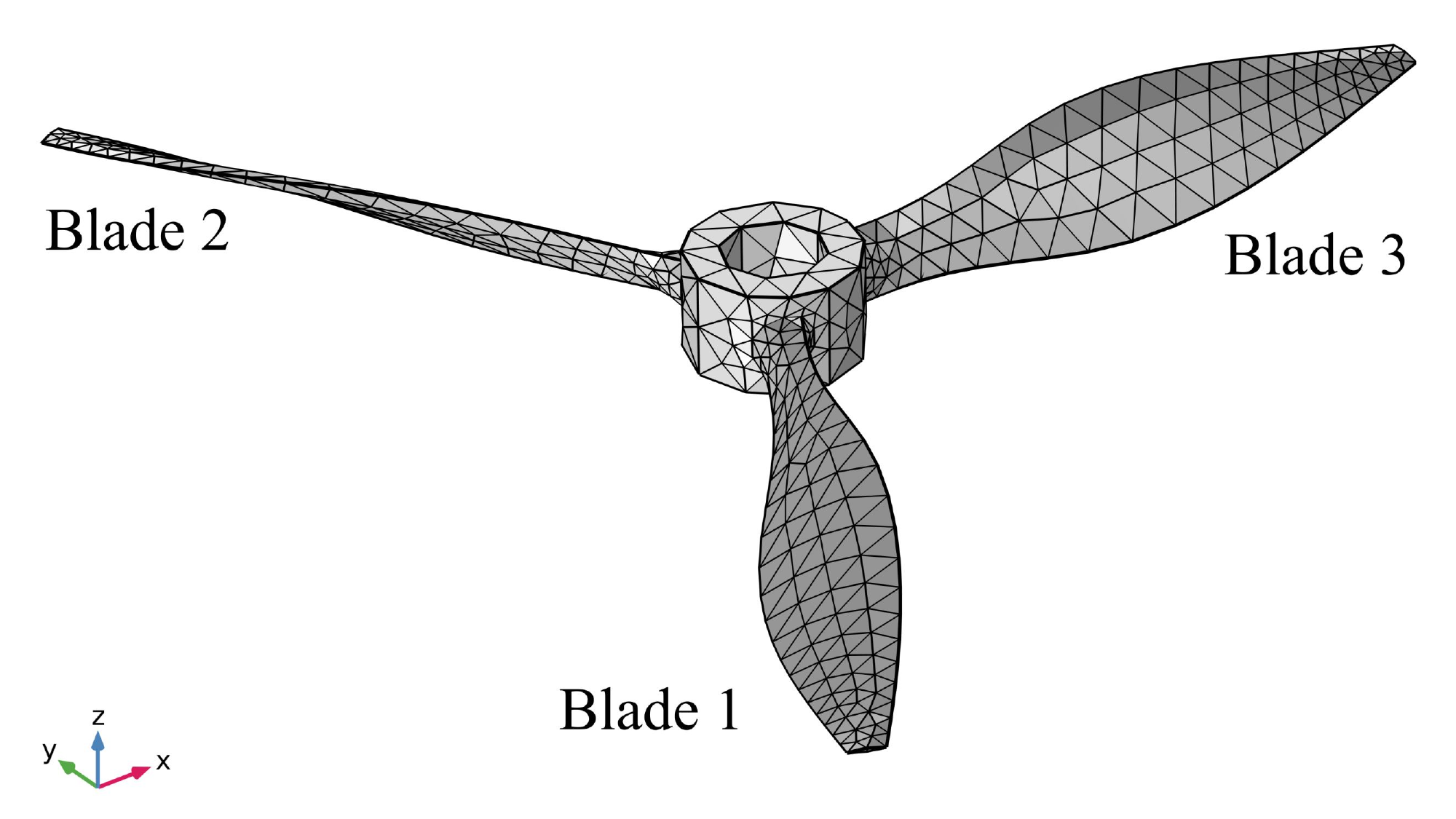}
	\caption{Geometry of the fan model and its finite element mesh. Three blades are labeled as the figure shows.}
	\label{fig-fan-model}
\end{figure}

\begin{figure*}
	\centering
	\includegraphics[width=170mm]{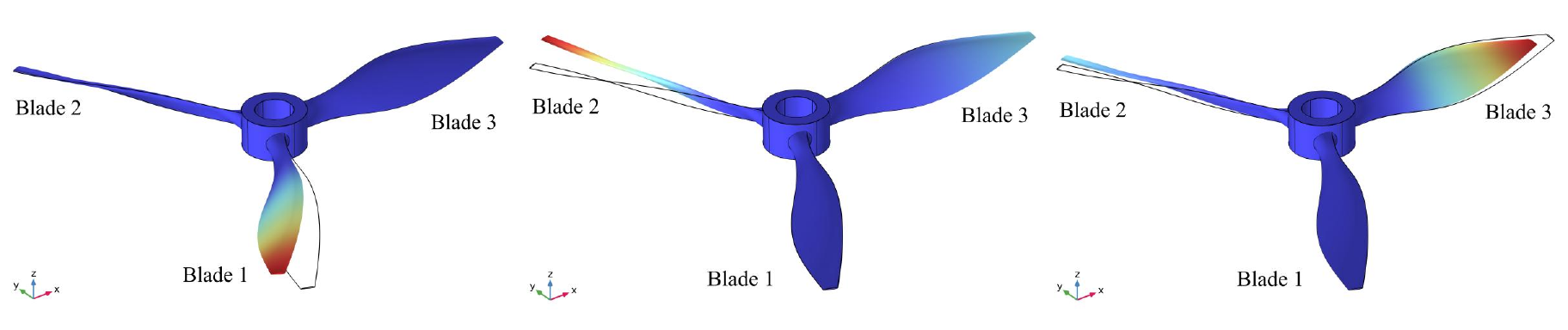}
	\caption{The first three mode shapes of the rotating fan. Left panel: 1st mode, middle panel: 2nd mode, right panel: 3rd mode.}
	\label{fig-fan-mode}
\end{figure*}

As our last example, we consider a rotating fan with three solid blades shown in \Cref{fig-fan-model}. We use this example to illustrate the effectiveness of our reduction method in rotating machines with complex geometry and multiple components. Here, each blade is approximately $2\,\text{m}$ in length, with a root width of $0.2\,\text{m}$. The blade cross-section adopts a modified NACA2410 airfoil, and twist angles are distributed along the spanwise direction. The central disk has an outer diameter of $0.5\,\text{m}$, an inner diameter of $0.3\,\text{m}$, and a thickness of $0.3\,\text{m}$, as illustrated in \Cref{fig-fan-model}. Fixed constraints are applied to the inner surface of the disk, and the rotational speed is set to $7.439\,\text{rad/s}$. The material properties adopted are identical to those used in \Cref{sec:beam}. Damping ratio of first mode $\xi_1$ is assigned as 0.001. The resulting finite element model has 11307 degrees of freedom. To facilitate reproducibility, the CAD geometry and finite element mesh used in this example have been made publicly available at~\cite{CAD-file}.

Owing to the cyclic symmetry of the structure, the frequencies of the first three modes for the whole structure are extremely close. In particular, we have $\omega_1=29.7545\,\text{rad/s}$, $\omega_2=29.7571\,\text{rad/s}$, and $\omega_3=29.7573\,\text{rad/s}$ for the undamped natural frequencies. Therefore, we have a near 1:1:1 internal resonance in this system. \Cref{fig-fan-mode} shows the mode shape of these modes. We observe that for each of these modes, one blade undergoes vibration similar to the first mode of a cantilever beam, while the other two blades are slightly deformed. In particular, we find that blades 1, 2, and 3 have relative large displacements for the 1st, 2nd, and 3rd modes, respectively. We note that these three mode shapes do not have perfect cyclic symmetry. In other words, rotating the mode vector of the first mode about the $z$-axis by $2\pi/3$ does not exactly yield the second mode.

\begin{figure*}
	\centering
	\includegraphics[width=170mm]{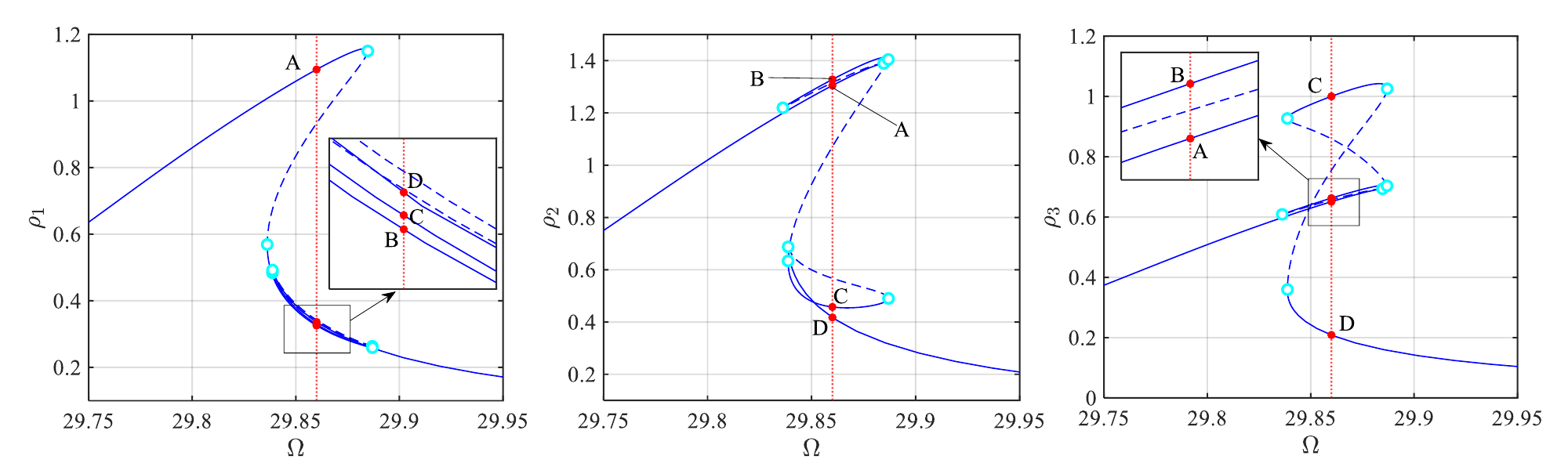}
	\caption{FRCs of SSM-based ROM for the simple fan model in reduced coordinates $\rho_1$ (left panel), $\rho_2$ (middle panel), and $\rho_3$ (right panel).}
	\label{fig-fan-rho}
\end{figure*}

\begin{figure*}
	\centering
	\includegraphics[width=170mm]{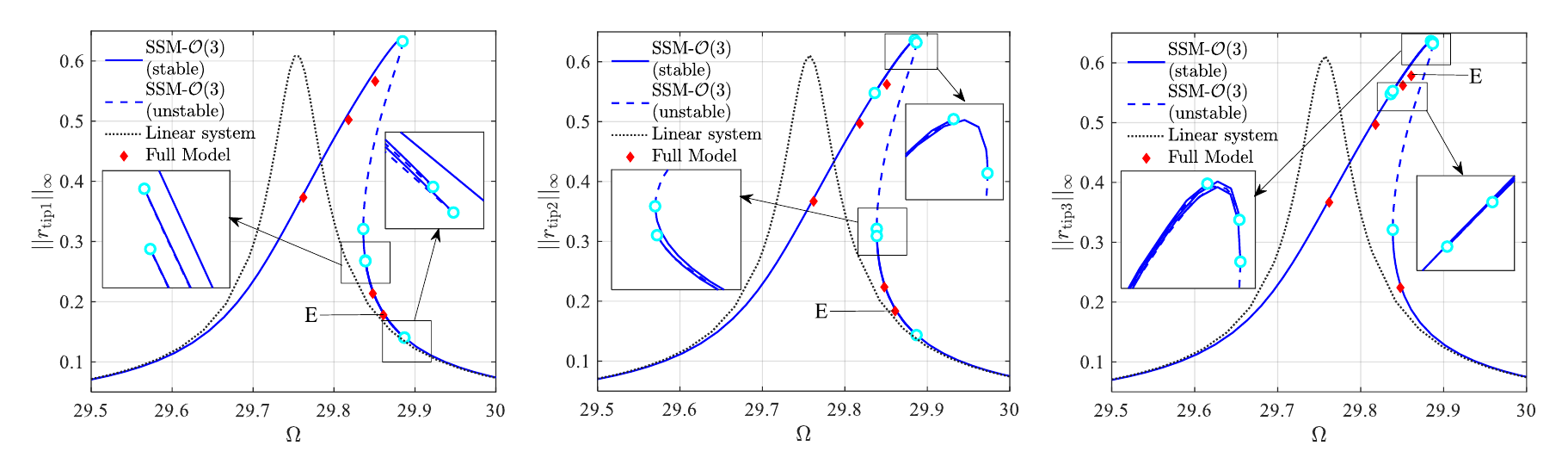}
	\caption{FRCs of the simple fan model for the total displacement at the tip of the three blades. Here, $r_{\mathrm{tip1}}$, $r_{\mathrm{tip2}}$, and $r_{\mathrm{tip3}}$ in the three panels denote the total displacement at the tip of the 1st, 2nd, and 3rd blades, respecitvely.}
	\label{fig-fan-frc}
\end{figure*}

Here, these internally resonant modes are then selected as the master subspace for model reduction, resulting in six-dimensional SSM-based ROM. We will use this SSM-based ROM to investigate the nonlinear response behavior of the system in the presence of internal resonance. Specifically, in-phase point loads $25\cos\Omega t$ are simultaneously applied at the tips of the three blades as external excitation to excite these three modes. 


We use the SSM-based ROM to extract FRCs of the system. Here, the reduced coordinates in polar form are given as $\bs{p}=(\rho_1e^{\mathrm{i}\theta_1},\rho_1e^{-\mathrm{i}\theta_1},\rho_2e^{\mathrm{i}\theta_2},\rho_2e^{-\mathrm{i}\theta_2},\rho_3e^{\mathrm{i}\theta_3},\rho_3e^{-\mathrm{i}\theta_3})$. We use the amplitudes $\rho_1$, $\rho_2$, and $\rho_3$ to present the FRCs in reduced coordinates. These FRCs are plotted in Fig.~\ref{fig-fan-rho} and enable us to explore the nonlinear interaction of the internally resonant modes. We observe six saddle-node (SN) bifurcation points denoted by cyan circles on each of the FRC. This complex bifurcation behavior is associated with the internal resonance. In particular, we find that for $\Omega\approx29.86$, seven periodic orbits coexist, among them there are four stable ones and three unstable ones. For the four coexisting stable periodic orbits marked as A-D in Fig.~\ref{fig-fan-rho}, the relative magnitudes for these reduced coordinates changes significantly. Specifically, we find that $\rho_1\approx\rho_2>\rho_3$ at point A, $\rho_2>\rho_1\approx\rho_3$ at point B, $\rho_3>\rho_1\approx\rho_2$ at point C, and $\rho_i$ is small for $i=1,2,3$ at point D.

We map these FRCs in reduced coordinates back to that in physical coordinates. Here, we present the FRCs for the total displacement at the tips of the three rotating blades. Notably, although the FRCs at the tips of the three blades appear highly similar in \Cref{fig-fan-frc}, the locations of SN points on these FRCs differ significantly. In particular, we find that the SN locations for the FRC at tip 1 are similar to those of FRC in $\rho_1$ (cf.~the upper panels of \Cref{fig-fan-frc} and \Cref{fig-fan-rho}), while the SN locations on the FRC at tip 2 or 3 is similar to that in $\rho_2$ or $\rho_3$. This indicates that the blades can have different amplitudes at the coexisting periodic orbits. Specifically, when one blade settles into the stable solution with a higher amplitude, the other blades may also reside in the high-amplitude stable solution or the low-amplitude stable solution.


Several representative excitation frequencies are chosen to calculate the responses of the full system for validation. The results show that the data obtained from the reduced-order system are in good agreement with those of the full model, as presented in \Cref{fig-fan-frc}. Among these representative points, the displacement at tip 3 is much larger than that of tips 1 and 2 at point E, which is qualitatively similar to point C in \Cref{fig-fan-rho}.
In this example, the total computation time for the SSM reduction is approximately $50\,\text{mins}$ due to increased master subspace dimensionality, while a single data point in \Cref{fig-fan-frc} would take about $5\,\text{hours}$ to obtain, which clearly demonstrates the high efficiency of SSM-based reduction.


\section{Conclusions}
\label{sec:conclusions} 
This study presents a non-intrusive model reduction framework for geometrically nonlinear rotating structures with Coriolis and centrifugal forces taken into consideration. This framework is based on the theory of spectral submanifold (SSM) and construct SSM-based reduced models (ROMs) associated with the nontrivial static equilibrium configuration induced by the centrifugal force. The effectiveness of this framework is demonstrated via a few representative examples, including a rotating beam, a twist blade, a rotor system, and an internally resonant fan.

SSM reduction effectively addresses the curse of high-dimensionality in the finite element model of rotating structures. Indeed, the SSM-based ROMs have only a few degrees of freedom and hence low-dimensional, yet make ensures the prediction accuracy for backbone and forced response curves. In addition, SSM reduction exhibits strong versatility and applicability in various rotating structures. It can accurately capture complex nonlinear behaviors such as geometric nonlinearity, centrifugal stiffening, and mode energy exchange induced by internal resonance. Further, the first two examples show that the Coriolis force can have a significant impact on the nonlinear vibration even though it has a limited impact on the natural frequency of rotating structures. This highlights the importance of Coriolis effects in nonlinear vibration analysis of rotating structures.

We have used COMSOL for finite element (FE) modeling in this study but note that the non-intrusive SSM reduction framework works well for other in-house or commercial FE code. Future extension of this study includes SSM reduction for more realistic applications such as blisk with tens of blades and the construction of parametric SSM-based ROMs that enable the prediction of responses under a range of rotating speeds. Another future research direction is the dynamics design such as vibration control and optimization of rotating structures via the SSM-based ROMs.

\section*{Acknowledgment}
In this work, Hejun Gao and Mingwu Li were supported by the National Natural Science Foundation of China (12572010) and the Pearl River Talent Recruitment Program (2023QN10H603).






\bibliographystyle{unsrt}
\bibliography{refs}

\end{document}